\documentclass[preprint,11pt]{elsarticle}

\usepackage{amssymb}
\usepackage{amsmath}
\usepackage{amsfonts}
\usepackage{amssymb}
\usepackage{indentfirst,latexsym,bm}
\usepackage{amsthm}
\usepackage{color,xcolor}
\usepackage[all]{xy}
\input{mathrsfs.sty}
\newtheorem{theorem}{Theorem}[section]
\newtheorem{lemma}[theorem]{Lemma}
\newtheorem{corollary}[theorem]{Corollary}
\newtheorem{proposition}[theorem]{Proposition}

\theoremstyle{definition}
\newtheorem{definition}{Definition}[section]

\theoremstyle{definition}

\theoremstyle{remark}
\newtheorem{remark}{Remark}[section]
\theoremstyle{question}

\theoremstyle{problem}

\numberwithin{equation}{section}

\journal{XXX}

\begin{document}

\begin{frontmatter}





\title{The block matrix representations for the quasi-projection pairs on Hilbert $C^*$-modules}
\author[shnu]{Xiaoyi Tian}
\ead{tianxytian@163.com}
\author[shnu]{Qingxiang Xu}
\ead{qingxiang\_xu@126.com}
\author[HS]{Chunhong Fu}
\ead{fchlixue@163.com}
\address[shnu]{Department of Mathematics, Shanghai Normal University, Shanghai 200234, PR China}
\address[HS]{Health School Attached to Shanghai University of Medicine $\&$ Health Sciences,   Shanghai 200237, PR China}

\begin{abstract}
A quasi-projection pair consists of two operators $P$ and $Q$ acting on a Hilbert $C^*$-module $H$, where $P$ is a projection and $Q$ is an idempotent satisfying
$Q^*=(2P-I)Q(2P-I)$, in which
$Q^*$ denotes the adjoint operator of $Q$, and $I$ is the identity operator on $H$. Such a pair is said to be harmonious if both $P(I-Q)$ and $(I-P)Q$ admit polar decompositions. The primary goal of this paper is to present the block matrix representations for a harmonious quasi-projection pair $(P,Q)$ on a Hilbert $C^*$-module, and additionally to derive new block matrix representations for the matched projection, the range projection, and the null space projection of $Q$. Several  applications of these newly obtained block matrix representations are also explored.
\end{abstract}


\end{frontmatter}




\tableofcontents

\section{Introduction}\label{sec:Intro}

Given  two (right) Hilbert modules $H$ and $K$ over a $C^*$-algebra, let $\mathcal{L}(H,K)$ be the set of the adjointable operators from $H$ to $K$, with the abbreviation $\mathcal{L}(H)$ whenever $H=K$, and the replacements of $\mathbb{B}(H,K)$ and $\mathbb{B}(H)$  in the special case where  $H$ and $K$ are Hilbert spaces.
The identity operator on $H$ is denoted by $I_H$,  or simply $I$ when no confusion arises.
By an idempotent $Q\in\mathcal{L}(H)$, we mean $Q^2=Q$.  If furthermore $Q$ is self-adjoint, then $Q$ is called a projection.
An idempotent that is not a projection is called a non-projection idempotent, while a projection $P$ is termed non-trivial if  $P\ne I$ and $P\ne 0$. It is notable that in some literatures, idempotents and projections are referred to be oblique/skew projections and orthogonal projections, respectively \cite{Afriat,BS,Ovchinnikov}.

It is known that up to unitary equivalence, every two projections $P$ and $Q$ on a Hilbert space have the  canonical $2\times 2$ block matrix representation \cite[Theorem~4.5]{Corach-Maestripieri}  and the Halmos' $6\times 6$ block matrix representation (see e.g.,\,\cite[Theorem~1.1]{BS} and \cite[Theorem~2]{Halmos}). These block matrix representations have numerous applications in operator theory. For instance, it is shown in \cite[Theorem~4.2]{Afriat}  that
\begin{equation}\label{formula for Q-001}Q=(I-P_{\mathcal{R}(Q)}P_{\mathcal{N}(Q)})^{-1}P_{\mathcal{R}(Q)} (I-P_{\mathcal{R}(Q)}P_{\mathcal{N}(Q)})
\end{equation} for every idempotent $Q\in\mathbb{B}(H)$,
where $P_{\mathcal{R}(Q)}$ and $P_{\mathcal{N}(Q)}$ denote the range projection and the null space projection of $Q$, respectively. Based on \eqref{formula for Q-001} and the Halmos' $6\times 6$ block matrix representations of $P_{\mathcal{R}(Q)}$ and $P_{\mathcal{N}(Q)}$, a useful block matrix representation of $Q$ can be derived; see \cite[Corollary~1.7]{BS} for the details.

Hilbert $C^*$-modules are the natural generalizations of Hilbert spaces by allowing the inner product to take values in a $C^*$-algebra rather than in the complex field $\mathbb{C}$ \cite{Lance,MT,Paschke}. A critical distinction from the Hilbert space is that a closed submodule of a Hilbert $C^*$-module  may fail to be orthogonally complemented.  This fundamental disparity leads to significant differences in operator properties: while every bounded linear operator on a Hilbert space admits the polar decomposition, the adjointable operators on Hilbert $C^*$-modules generally do not (see e.g.,\,\cite[Example~3.15]{Liu-Luo-Xu}). Consequently, the canonical $2\times 2$ block matrix representation and the Halmos' $6\times 6$ block matrix representation, which rely on the polar decompositions of operators, become inapplicable for a general pair of projections $P$ and $Q$ on a Hilbert $C^*$-module. Actually, it is proved in \cite[Theorems~2.14 and 3.3]{Xu-Yan} that these block matrix representations are applicable if and only if the pair $(P, Q)$ of projections satisfies the harmonious condition defined in \cite[Definition~4.1]{Luo-Moslehian-Xu} and furthermore characterized in \cite[Lemmas~2.12 and 2.13]{Xu-Yan} (see also \cite[Theorem~2.4]{FXY}).

Let $H$ be a Hilbert $C^*$-module. Motivated by the investigation of the operator distances from projections in $\mathcal{L}(H)$ to a given idempotent $Q\in\mathcal{L}(H)$,  a specific projection $m(Q)$,  called the matched projection of $Q$,  is introduced  in \cite{TXF}. Inspired by the relationship between an idempotent $Q$ and its matched projection $m(Q)$,
another concept called the quasi-projection pair $(P,Q)$ is also introduced in \cite{TXF}, where $P$ is a projection, while $Q$ is an idempotent  satisfying
$$Q^*=(2P-I)Q(2P-I).$$
According to
\cite[Definition~3.2]{TXF-03}, a quasi-projection pair $(P,Q)$ on a Hilbert $C^*$-module $H$ is said to be semi-harmonious if both
$P(I-Q)$ and $(I-P)Q$  have the polar decompositions. Furthermore, $(P,Q)$ is said to be harmonious if  both
$\overline{\mathcal{R}\big(PQ(I-P)\big)}$ and $\overline{\mathcal{R}\big((I-P)QP\big)}$  are  orthogonally complemented in $H$
\cite[Definition~4.1]{TXF-03}. Some detailed characterizations of the semi-harmonious and harmonious quasi-projection pairs on Hilbert $C^*$-modules are given in \cite{TXF-03}, establishing the technical groundwork for studying their block matrix representations.

The main purpose of this paper is to present the block matrix representations for a harmonious quasi-projection pair $(P,Q)$ on a Hilbert $C^*$-module, and additionally to derive new block matrix representations for the associated projections $m(Q), P_{\mathcal{R}(Q)}$, and $P_{\mathcal{N}(Q)}$, since, as shown in \eqref{equ:final exp of P}, \eqref{range projection Q} and \eqref{equ:relation of Q and P NQ}, these projections can be expressed via operators in the $C^*$-algebra generated by the idempotent $Q$ and the identity operator $I$. Let $(P,Q)$ be a harmonious quasi-projection pair on a Hilbert $C^*$-module $H$. Up to unitary equivalence, the canonical $2\times 2$ block matrix representation  and the Halmos-like  $6\times 6$  block matrix representation of $(P,Q)$ are established in Theorems~\ref{thm:2 by 2 canonical representation} and \ref{thm:the 6x6 representation of p and q}, respectively.
 Based on the block matrix representation of $Q$ obtained in Theorem~\ref{thm:2 by 2 canonical representation}, new $2\times 2$ block matrix representations of the above three projections are derived in
Corollaries~\ref{cor:derived rep for range Q} and \ref{cor:derived rep for kernel Q}, and Theorem~\ref{thm:matrix rep for m Q when that Q is given}, respectively. Notably, new $6\times 6$ block matrix representations of these projections can also be derived by using  the block matrix representation of $Q$ from Theorem~\ref{thm:the 6x6 representation of p and q}. It is remarkable that $\big(m(Q),Q\big)$ is always a quasi-projection pair for every idempotent $Q\in\mathcal{L}(H)$. So, if $\big(m(Q),Q\big)$ is harmonious (which occurs when $H$ is a Hilbert space), then as shown in Theorems~\ref{thm:representation of Q realate to mQ}, \ref{thm:new representation for a general idempotent} and \ref{thm:new 6order rep for range and null}, new $2\times 2$ and $4\times 4$ block matrix representations of $Q$, $P_{\mathcal{R}(Q)}$ and $P_{\mathcal{N}(Q)}$ can be derived.

This paper further presents several applications. For every idempo\-tent $Q\in\mathcal{L}(H)$, a projection $s(Q)$, termed the supplementary projection of $Q$, is introduced and clarified in Section~\ref{subsec:supplementary projection}. As an application, Theorem~\ref{thm:new formula for q bansed on mQ and sQ} provides a new formula expressing $Q$ via  the projections $m(Q)$ and $s(Q)$, rather than the projections $P_{\mathcal{R}(Q)}$ and $P_{\mathcal{N}(Q)}$. Using Theorem~\ref{thm:new representation for a general idempotent}, a new proof of the canonical form for quadratic operators is presented (see Theorem~\ref{thm:quadratic form}). The Krein-Krasnoselskii-Milman equality is a well-known norm identity, which is particularly useful in dealing with the norms associated with projection pairs. It is interesting to investigate the validity of Krein-Krasnoselskii-Milman equality in the context of quasi-projection pairs. Based on the block matrix representations of the associated operators, Theorems~\ref{thm:modified KKM formula} and \ref{thm:KKME for matched pair} show that this equality does not hold in general for quasi-projection pairs.

Unless otherwise specified, throughout the remainder of this paper, $\mathbb{N}$ is the set of positive integers, $\mathbb{R}$ is the real line, $\mathbb{C}$ is the complex field, $\mathfrak{A}$ is a $C^*$-algebra, $\mathfrak{A}_{\mbox{sa}}$ is the set of self-adjoint elements in  $\mathfrak{A}$,  $E,H$ and $K$ are Hilbert modules over $\mathfrak{A}$, $\mathcal{L}(H)_{\mbox{sa}}$ (resp.\,$\mathcal{L}(H)_+$) is the set of self-adjoint (resp.\,positive) elements in the $C^*$-algebra $\mathcal{L}(H)$. For each $A\in\mathcal{L}(H)$, the notation $A\ge 0$ means that $A$ is a positive operator on $H$.

It is known that every adjointable operator $A\in\mathcal{L}(H,K)$ is a bounded linear operator and is also $\mathfrak{A}$-linear in the sense that
\begin{equation*}\label{equ:keep module operator}A(xa)=A(x)a,\quad\forall\,x\in H, a\in\mathfrak{A}.
\end{equation*}
Let $\mathcal{R}(A)$, $\mathcal{N}(A)$, and $A|_M$ denote the rang of $A$, the null space of $A$, and the restriction of $A$ to a subset $M$ of $H$, respectively.
When $A\in\mathcal{L}(H)$ or $A\in \mathfrak{A}$,
   let $\sigma(A)$ denote the spectrum of $A$ with respect to $\mathcal{L}(H)$ or $\mathfrak{A}$.
A closed submodule $M$ of  $H$ is said to be
orthogonally complemented  in $H$ if $H=M+ M^\bot$, where
$$M^\bot=\big\{x\in H:\langle x,y\rangle=0\ \mbox{for every}\ y\in
M\big\}.$$
In this case, the projection from $H$ onto $M$ is denoted by $P_M$.
Let $v^T$ denote the transpose of a vector $v$ in a partitioned linear space. Define
$$H\oplus K=\big\{(x, y)^T :x\in H,  y\in K\big\},$$ which is also a Hilbert $\mathfrak{A}$-module with the  $\mathfrak{A}$-valued inner product given by
$$\left<(x_1, y_1)^T, (x_2, y_2)^T\right>=\big<x_1,x_2\big>+\big<y_1,
y_2\big>$$ for all  $x_1,x_2\in H$ and $y_1,y_2\in K$.

The rest of this paper is organized  as follows.  Section~\ref{sec:some lemmas} introduces four frequently used functions--$f,g,h$, and $\ell$, defined on the subset $(0,1)^c$ of the real line, and provides basic knowledge about the polar decomposition, the quasi-projection pair, and the matched projection. Sections~\ref{sec:representations for qpp} and \ref{sec:66 representations for qpp} are devoted to the canonical $2\times 2$ block matrix representation, and the Halmos-like block matrix representation  for a harmonious quasi-projection pair, respectively. Section~\ref{sec:applications} explores the applications of these newly obtained block matrix representations.

\section{Some preliminaries}\label{sec:some lemmas}

The functional calculus for elements in a commutative $C^*$-algebra will be used frequently in this paper, so it is helpful to specify the functions employed in the sequel. For this purpose, we define
$$(0,1)^c=(-\infty,0]\cup [1,+\infty),$$ and let $f,g,h$ and $\ell$ be the real-valued functions defined on $(0,1)^c$ by
\begin{align}&\label{defn of f-ist time}f(t)=\begin{cases}1, &t\in [1,+\infty),\\
-1,&t\in (-\infty, 0],
\end{cases}\\
&\label{equ: the function of g,h,ell}g(t)=\sqrt{|t|},\quad  h(t)=-f(t)\sqrt{|t-1|},\quad \ell(t)=\sqrt{t^2-t}.
\end{align}
Clearly, all these functions are continuous and satisfy the following identities for every $t\in (0,1)^c$:
\begin{align*}&f^2(t)\equiv 1,\quad h^2(t)+g^2(t)=|2t-1|, \\
  &h(t)t+g(t)\ell(t)\equiv 0, \quad g^2(t)-h^2(t)=f(t),\\
&h(t)\ell(t)+g(t)(t-1)\equiv 0.
\end{align*}
Consequently, if  $A$ is a self-adjoint element in a unital $C^*$-algebra $\mathfrak{A}$ such that
$\sigma(A)\subseteq (0,1)^c$, then
\begin{equation}\label{defn of operator ell A}
\ell(A)=(A^2-A)^\frac{1}{2},
\end{equation}
and the following operator identities hold:
\begin{align}
\label{equ:newtemp-003}&f^2(A)=I_{\mathfrak{A}},\quad h^2(A)+g^2(A)=|2A-I_{\mathfrak{A}}|,\\
&\label{equ:the properties of f and g}
  h(A)A+g(A)\ell(A)=0, \quad g^2(A)-h^2(A)=f(A),\\
&\label{equ:temp-002}
  h(A)\ell(A)+g(A)(A-I_{\mathfrak{A}})=0.
\end{align}

To establish the main results of this paper, we first recall the following known lemma.

\begin{lemma}\label{lem:keep orthogonal complementity}{\rm \cite[Lemma~2.1]{Xu-Yan}} Let $M$ and $N$ be closed submodules of $H$ such that $N\subseteq M$ and $M$ is orthogonally complemented in $H$. Then the following statements are equivalent:
\begin{enumerate}
\item[{\rm (i)}] $N$ is orthogonally complemented in $H$;
\item[{\rm (ii)}] $N$ is orthogonally complemented in $M$.
\end{enumerate}
\end{lemma}

For every $T\in\mathcal{L}(H,K)$, let $|T|$ denote the (positive) square root of $T^*T$. Then $|T^*|=(TT^*)^\frac12$. According to  \cite[Lemma~3.9]{Liu-Luo-Xu}, there exists at most one partial isometry
$U\in\mathcal{L}(H,K)$  satisfying
\begin{equation}\label{equ:two conditions of polar decomposition}T=U|T|\quad \mbox{with}\quad U^*U=P_{\overline{\mathcal{R}(T^*)}}.\end{equation}
This representation is referred to as the polar decomposition of $T$  \cite[Definition~3.10]{Liu-Luo-Xu}. However,  as illustrated in  \cite[Example~3.15]{Liu-Luo-Xu}, an adjointable operator on a Hilbert $C^*$-module may not always admit a polar decomposition. The following lemma provides a criterion for its existence.

\begin{lemma}\label{lem:polar decomposition of T star}{\rm \cite[Lemma~3.6 and Theorem~3.8]{Liu-Luo-Xu}} For every $T\in\mathcal{L}(H,K)$, the following statements are  equivalent:
\begin{enumerate}
\item[{\rm (i)}] $T$ has the polar decomposition;
\item[{\rm (ii)}]$T^*$ has the polar decomposition;
\item[{\rm (iii)}]$\overline{\mathcal{R}(T)}$ and $\overline{\mathcal{R}(T^*)}$ are orthogonally complemented in $K$ and $H$, respectively.
\end{enumerate}
If any of the conditions {\rm (i)--(iii)} is satisfied, and the polar decomposition of $T$ is given by \eqref{equ:two conditions of polar decomposition}, then the polar decomposition of $T^*$ is given  by
\begin{equation*}\label{equ:the polar decomposition of T star-pre stage}T^*=U^*|T^*|\quad \mbox{with}\quad UU^*=P_{\overline{\mathcal{R}(T)}}.
\end{equation*}
\end{lemma}

\begin{remark}
Suppose $T\in\mathcal{L}(H,K)$ has the representation \eqref{equ:two conditions of polar decomposition}. For simplicity, we hereafter say that $T$ admits polar decomposition $T=U|T|$.
\end{remark}

\begin{lemma}\label{lem:rang characterization-1}{\rm \cite[Proposition~2.7]{Liu-Luo-Xu}} Let $A\in\mathcal{L}(H,K)$ and $B,C\in\mathcal{L}(E,H)$ satisfy  $\overline{\mathcal{R}(B)}=\overline{\mathcal{R}(C)}$. Then $\overline{\mathcal{R}(AB)}=\overline{\mathcal{R}(AC)}$.
\end{lemma}

\begin{lemma}\label{lem:Range Closure of T alpha and T}{\rm (\cite[Proposition~2.9]{Liu-Luo-Xu} and \cite[Lemma 2.2]{Vosough-Moslehian-Xu})} Let $T\in \mathcal{L}(H)_+$. Then for every $\alpha>0$, $\overline{\mathcal{R}(T^{\alpha})}=\overline{\mathcal{R}(T)}$ and $\mathcal{N}(T^{\alpha})=\mathcal{N}(T)$.
\end{lemma}

\begin{lemma}\label{lem:Range closure of TT and T} {\rm\cite[Proposition 3.7]{Lance}}
For every $T\in\mathcal{L}(H,K)$, $\overline {\mathcal{R}(T^*T)}=\overline{ \mathcal{R}(T^*)}$ and $\overline {\mathcal{R}(TT^*)}=\overline{ \mathcal{R}(T)}$.
\end{lemma}

Next, we introduce the fundamental concepts of the quasi-projection pair and the matched projection.

\begin{definition}\label{defn:quasi-projection pair}\rm{ \cite[Definition~2.2]{TXF}} An ordered pair $(P,Q)$  is called a quasi-projection pair on $H$ if $P\in\mathcal{L}(H)$ is a projection and $Q\in\mathcal{L}(H)$ is an idempotent such that
\begin{align}\label{conditions 1 for qpp}&PQ^*P=PQP, \quad PQ^*(I-P)=-PQ(I-P),\\
\label{conditions 2 for qpp}&(I-P)Q^*(I-P)=(I-P)Q(I-P).
\end{align}
\end{definition}

We present two useful characterizations of quasi-projection pairs as follows.

 \begin{lemma}\label{thm:four equivalences} \rm{ \cite[Theorem~2.4]{TXF}} Suppose that $P\in\mathcal{L}(H)$ is a projection and $Q\in\mathcal{L}(H)$ is an idempotent.  If one element in
  \begin{equation*}\label{Sigma eight elements}\Sigma:=\big\{(A,B): A\in \{P, I-P\}, B\in \{Q,Q^*, I-Q,I-Q^*\}\big\}\end{equation*}
  is a  quasi-projection pair, then all the remaining elements in $\Sigma$ are the quasi-projection pairs.
 \end{lemma}

\begin{lemma}\label{thm:short description of qpp}{\rm \cite[Theorem~2.5]{TXF}} Suppose that $P\in \mathcal{L}(H)$ is a projection and $Q\in\mathcal{L}(H)$ is an idempotent. Then the following statements are equivalent:
\begin{itemize}
  \item [\rm{(i)}] $(P,Q)$ is a quasi-projection pair;
  \item [\rm{(ii)}]$Q^*=(2P-I)Q(2P-I)$;
   \item [\rm{(iii)}]$|Q^*|=(2P-I)|Q|(2P-I)$.
 \end{itemize}
\end{lemma}

\begin{remark}\label{rem:two projections qpp}Let $P,Q\in\mathcal{L}(H)$ be two projections. From the equivalence of (i) and (ii) in the preceding lemma, it follows that $(P,Q)$ is a quasi-projection pair if and only if $PQ=QP$, which is equivalent to $PQ$ being a projection.
\end{remark}

Before proceeding, we introduce the following notations for operators and associated submodules.

\begin{definition}\label{defn of 6 operators and modules} Given two idempotents $P,Q\in \mathcal{L}(H)$,  define the operators:
\begin{align}&T_1=P(I-Q),  \quad    T_2=(I-P)Q,\nonumber\\
\label{equ:new defn of T3 and T4}&T_3=PQ(I-P),\quad T_4=(I-P)QP,
 \end{align}
and define the submodules:
\begin{align}\label{eqn:defn of H1 and H4}&H_1=\mathcal{R}(P)\cap\mathcal{R}(Q), \quad H_4=\mathcal{N}(P)\cap\mathcal{N}(Q),\\
 \label{eqn:defn of H2 and H3}& H_2=\mathcal{R}(P)\cap\mathcal{N}(Q),  \quad H_3=\mathcal{N}(P)\cap\mathcal{R}(Q), \\
 \label{eqn:defn of H5 and H6}&H_5=\overline{\mathcal{R}(T_3)}, \quad  H_6=\overline{\mathcal{R}(T_4)}.
\end{align}
\end{definition}

The following lemma establishes the equivalence of the conditions $H_5=0$ and $H_6=0$.
\begin{lemma}\label{lem:H5 is 0}Let $(P,Q)$ be a quasi-projection pair. Define  $H_5$ and $H_6$ as in \eqref{eqn:defn of H5 and H6}. Then
the following statements are equivalent:
\begin{enumerate}
\item[{\rm (i)}] $Q=Q^*$;
\item[{\rm (ii)}] $H_5=\{0\}$;
\item[{\rm (iii)}]$H_6=\{0\}$.
\end{enumerate}
\end{lemma}
\begin{proof}(i)$\Longrightarrow$(ii). By Remark~\ref{rem:two projections qpp}, we have $PQ=QP$. In this case, the operator $T_3$ defined by
\eqref{equ:new defn of T3 and T4} is zero, which leads by \eqref{eqn:defn of H5 and H6} to $H_5=\{0\}$.

(ii)$\Longrightarrow$(iii). From the definitions of $H_5$ and $T_3$, we have $PQ=PQP$.  Taking adjoints of both sides yields
$$(2P-I)Q(2P-I)P=PQP.$$
Simplifying this expression, we obtain $(P-I)QP=0$, that is, $-T_4=0$. Hence, by \eqref{eqn:defn of H5 and H6}, we conclude that $H_6=\{0\}$.

(iii)$\Longrightarrow$(i). By assumption, we have $(I-P)QP=0$. Taking adjoints gives
$$0=P(2P-I)Q(2P-I)(I-P)=-PQ(I-P).$$
It follows that
$$QP=PQP=PQ.$$
Hence,
$$Q^*=(2P-I)Q(2P-I)=Q(2P-I)^2=Q.$$
This completes the proof.
\end{proof}

Regarding the operator ranges associated with a quasi-projection pair, we present the following lemma.
 \begin{lemma}\label{lem:another form of H1 and H4}{\rm \cite[Lemma~3.3]{TXF-03}} Suppose that $(P,Q)$ is a quasi-projection pair on $H$. Let $H_i$ $(1\le i\le 4)$ be defined by  \eqref{eqn:defn of H1 and H4}--\eqref{eqn:defn of H2 and H3}. Then
\begin{align}
&\label{alternative exps of H1 and H4}H_1=\mathcal{R}(P)\cap \mathcal{R}(Q^*),\quad H_4=\mathcal{N}(P)\cap \mathcal{N}(Q^*),\\
&\label{alternative exps of H2 and H3}H_2=\mathcal{R}(P)\cap \mathcal{N}(Q^*),\quad H_3=\mathcal{N}(P)\cap \mathcal{R}(Q^*).
\end{align}
\end{lemma}

Recall that the Moore-Penrose inverse $T^\dag$ of an operator $T\in \mathcal{L}(H,K)$  is the
unique element $X\in \mathcal{L}(K,H)$ satisfying
\begin{equation*} \label{equ:m-p inverse} TXT=T,\quad XTX=X,\quad (TX)^*=TX, \quad (XT)^*=XT.\end{equation*}
If such an operator  $T^\dag$ exists, then  $T$ is said to be Moore-Penrose invertible. For every $T\in\mathcal{L}(H,K)$, it follows from \cite[Theorem~2.2]{XS} that $T$ is Moore-Penrose invertible if and only if $\mathcal{R}(T)$ is closed in $K$.

Now, we describe the concept of the matched projection and present two of its key properties.

\begin{definition}\cite[Definition~3.2 and Theorem~3.4]{TXF}\label{def of mq} For every idempotent $Q\in \mathcal{L}(H)$, its matched projection is defined as
\begin{equation}\label{equ:final exp of P}m(Q)=\frac12\big(|Q^*|+Q^*\big)|Q^*|^\dag\big(|Q^*|+I\big)^{-1}\big(|Q^*|+Q\big),
\end{equation}
where $|Q^*|^\dag$ denotes the Moore-Penrose inverse of $|Q^*|$, which can be expressed  as
\begin{equation}\label{formula for the MP inverse-01}|Q^*|^\dag=\left(P_{\mathcal{R}(Q)}P_{\mathcal{R}(Q^*)}P_{\mathcal{R}(Q)}\right)^\frac12.
\end{equation}
\end{definition}

In the remainder  of this paper, the notations $m(Q)$ and $\big(m(Q),Q\big)$ denote the matched projection and
 the matched pair of an idempotent $Q$, respectively. It is known  that for every idempotent $Q$,  the pair $\big(m(Q),Q\big)$ is a quasi-projection pair \cite[Theorem~3.1]{TXF}.

\begin{remark}\label{remark of Q 02} Suppose that $Q$ is a non-projection idempotent. Up to unitary equivalence, $Q$ can be expressed as
\begin{equation}\label{defn of widetile Q111}Q=\left(
                                                       \begin{array}{cc}
                                                         I_{\mathcal{R}(Q)} & A \\
                                                         0 & 0 \\
                                                       \end{array}
\right)\in \mathcal{L}\big(\mathcal{R}(Q)\oplus\mathcal{N}(Q^*)\big)\end{equation}
for some $A\in \mathcal{L}\big(\mathcal{N}(Q^*), \mathcal{R}(Q)\big)$. Based on \eqref{equ:final exp of P}, it can be shown (see e.g.,\,\cite[Theorem~3.1]{TXF}) that
\begin{equation}\label{defn of widetilde Q wrt A B}
m(Q)=\frac12 \left(\begin{array}{cc}
        \big(B+I_{\mathcal{R}(Q)}\big)B^{-1} &  B^{-1}A\\
        A^*B^{-1} &  A^*\big[B\big(B+I_{\mathcal{R}(Q)}\big)\big]^{-1}A \\
      \end{array}\right),
\end{equation}
where $B=(AA^*+I_{\mathcal{R}(Q)})^{\frac{1}{2}}$.

\end{remark}

\begin{lemma}\label{thm:matched projection for I-Q}{\rm \cite[Theorems 3.7 and 3.14]{TXF}} Let $Q\in \mathcal{L}(H)$ be an idempotent. Then $$m(Q^*)=m(Q), \quad  m(I-Q)=I-m(Q).$$
\end{lemma}

\begin{lemma}\label{lem:matched pair H2 and H3 zero}{\rm \cite[Corollary~3.15]{TXF}} Let $Q\in \mathcal{L}(H)$ be an idempotent. Then
\begin{equation*}\label{two intersections are zero}\mathcal{R}\big(m(Q)\big)\cap \mathcal{N}(Q)=\{0\},\quad \mathcal{N}\big(m(Q)\big)\cap \mathcal{R}(Q)=\{0\}.\end{equation*}
\end{lemma}

\section{The $2\times 2 $ block matrix representations for harmonious  quasi-projection pairs}\label{sec:representations for qpp}

  This section presents the canonical $2\times 2$ block matrix representation for a general harmonious quasi-projection pair $(P,Q)$. This presentation enables us to derive new $2\times 2$ block matrix representations for the associated projections  $P_{\mathcal{R}(Q)}$,  $P_{\mathcal{N}(Q)}$, and $m(Q)$. Furthermore, in the special case where the matched pair  $\big(m(Q),Q\big)$ is harmonious, we establish a new $2\times 2$ block matrix representation of $Q$.

\begin{definition}\label{defn of harmonious quasi-projection pair}{\rm\cite[Definition~4.1]{TXF-03}} A quasi-projection pair $(P,Q)$ on $H$ is said to be harmonious if the submodules $H_5$ and $H_6$ defined by \eqref{eqn:defn of H5 and H6}
are both orthogonally complemented in $H$.
\end{definition}

 Given a non-trivial projection $P\in\mathcal{L}(H)$,  let $U_P: H\to \mathcal{R}(P)\oplus \mathcal{N}(P)$ be the unitary operator defined  by
\begin{equation}\label{equ:unitary operator induced by P}U_Ph=\binom{Ph}{(I-P)h},\quad \forall\, h\in H.
\end{equation}
Clearly, for every $h_1\in\mathcal{R}(P)$ and $h_2\in\mathcal{N}(P)$,
\begin{equation*}\label{expression of inverse of U P}U_P^{*}\binom{h_1}{h_2}=h_1+h_2.\end{equation*}
Therefore, for every $T\in \mathcal{L}(H)$,  we have  (see e.g.,\,\cite[Section~2.2]{Xu-Wei-Gu})
 \begin{equation}\label{equ:block matrix T}
U_PTU_P^{*}=\left(
            \begin{array}{cc}
              PTP|_{\mathcal{R}(P)} & PT(I-P)|_{\mathcal{N}(P)} \\
              (I-P)TP|_{\mathcal{R}(P)} &(I-P)T(I-P)|_{\mathcal{N}(P)}\\
            \end{array}
          \right).
\end{equation}
In particular,
 \begin{equation}\label{equ:block matrix P}
U_PPU_P^{*}=\left(
            \begin{array}{cc}
              I_{\mathcal{R}(P)} & 0 \\
              0 &0\\
            \end{array}
          \right),\quad U_P(I-P)U_P^{*}=\left(
            \begin{array}{cc}
              0 & 0 \\
              0 &I_{\mathcal{N}(P)}\\
            \end{array}
          \right).
\end{equation}
For any $Q\in\mathcal{L}(H)$, by \eqref{equ:block matrix T},  \eqref{conditions 1 for qpp} and
\eqref{conditions 2 for qpp},  we conclude that $(P,Q)$ is
a quasi-projection pair if and only if
\begin{equation}\label{equ:block matrix Q}U_PQU_P^{*}=\left(
              \begin{array}{cc}
                A & -Q_{21}^* \\
                Q_{21} & Q_{22} \\
              \end{array}
            \right),
\end{equation}
where
$A\in \mathcal{L}\big(\mathcal{R}(P)\big)_{\mbox{sa}}$, $Q_{22}\in\mathcal{L}\big(\mathcal{N}(P)\big)_{\mbox{sa}}$ and $Q_{21}\in\mathcal{L}\big(\mathcal{R}(P),\mathcal{N}(P)\big)$ are such that $U_PQU_P^{*}$ is an idempotent.

Suppose that $P\in\mathcal{L}(H)$ is a non-trivial projection  and $Q\in\mathcal{L}(H)$  is an idempotent. Let $T_3$, $T_4$, $H_5$ and $H_6$  be defined by \eqref{equ:new defn of T3 and T4} and \eqref{eqn:defn of H5 and H6}, respectively. Since  $Q_{21}=T_4|_{\mathcal{R}(P)}$, by Lemma~\ref{lem:keep orthogonal complementity}  $\overline{\mathcal{R}(Q_{21})}$ is orthogonally complemented in $\mathcal{N}(P)$ if and only if $H_6$ is orthogonally complemented in $H$. Similarly,  $\overline{\mathcal{R}(Q_{21}^*)}$ is orthogonally complemented in $\mathcal{R}(P)$ if and only if
$H_5$ is orthogonally complemented in $H$. Combining the above observations with  Definition~\ref{defn of harmonious quasi-projection pair} and Lemma~\ref{lem:polar decomposition of T star}, we obtain
\begin{equation}\label{equ:existence of polar decomposition-01}\mbox{$(P,Q)$ is harmonious}\Longleftrightarrow \mbox{$Q_{21}$ has the polar decomposition.}\end{equation}

\begin{definition}\label{defn:canonical rep}Let $P\in\mathcal{L}(H)$ be a non-trivial projection. An operator $Q\in\mathcal{L}(H)$ is said to have the canonical $2\times 2$ block matrix representation induced by $P$ if
\begin{equation}\label{equ:UpQUp*}
U_P QU_P ^{*}=\left(
                \begin{array}{cc}
                  A & -\ell(A)U^* \\
                  U\ell(A) & U(I_{\mathcal{R}(P)}-A)U^*+Q_0 \\
                \end{array}
              \right),
\end{equation}
where $U_P$ is the unitary operator defined by \eqref{equ:unitary operator induced by P}, $A\in\mathcal{L}\big(\mathcal{R}(P)\big)_{\mbox{sa}}$ satisfies $A^2-A\ge 0$, $\ell(A)$ is defined by \eqref{defn of operator ell A},
$U\in \mathcal{L}\big(\mathcal{R}(P),\mathcal{N}(P)\big)$ is a partial isometry, and $Q_0\in \mathcal{L}\big(\mathcal{N}(P)\big)$ is a projection such that
\begin{equation}\label{equ:RU*=R Q112-Q11 U*Q=0}
\mathcal{R}(U^*)=\overline{\mathcal{R}(A^2-A)} \quad\mbox{and}\quad  U^*Q_0=0.
\end{equation}
\end{definition}

The first main result of this section is stated as follows.
\begin{theorem}\label{thm:2 by 2 canonical representation} Let $P$ be a non-trivial projection. Then for every $Q\in\mathcal{L}(H)$, the operator $Q$ admits the canonical $2\times 2$ block matrix representation induced by $P$ if and only if $(P,Q)$ is a harmonious quasi-projection pair.
\end{theorem}
\begin{proof} Let $T_3$, $T_4$, $H_5$, and $H_6$  be defined by \eqref{equ:new defn of T3 and T4} and \eqref{eqn:defn of H5 and H6}, respectively. For simplicity, we denote  $I_{\mathcal{R}(P)}=I_1$ and $I_{\mathcal{N}(P)}=I_2$.

$\Longleftarrow$.  Let  $U_P QU_P ^{*}$ be given by \eqref{equ:UpQUp*} such that the operators $A$, $U$, and $Q_0$ satisfy the conditions stated in Definition~\ref{defn:canonical rep}.
By Lemma~\ref{lem:Range Closure of T alpha and T} and the first equation in \eqref{equ:RU*=R Q112-Q11 U*Q=0},   we have
\begin{equation}\label{equ:many ranges are equal-02}\overline{\mathcal{R}[\ell(A)]}=\overline{\mathcal{R}[\ell^2(A)]}=\overline{\mathcal{R}(A^2-A)}=\mathcal{R}(U^*)=\mathcal{R}(U^*U).\end{equation}
Utilizing Lemma~\ref{lem:Range closure of TT and T} gives   \begin{equation*}\overline{\mathcal{R}\big[(I_1-A)^2\big]}=\overline{\mathcal{R}\big[(I_1-A)(I_1-A)^*\big]}=\overline{\mathcal{R}(I_1-A)}.\end{equation*}
Combining the above equations with \eqref{equ:many ranges are equal-02} and Lemma~\ref{lem:rang characterization-1}, we obtain
\begin{align}
  \overline{\mathcal{R}\big[(I_1-A)U^*\big]} & =\overline{\mathcal{R}\big[(I_1-A)(A^2-A)\big]}=\overline{\mathcal{R}\big[A(I_1-A)^2\big]}  \nonumber\\
 \label{range contained in that of T-01} & =\overline{\mathcal{R}\big[A(I_1-A)\big]}=\overline{\mathcal{R}[-\ell^2(A)]}=\overline{\mathcal{R}[\ell(A)]}.
\end{align}
From \eqref{equ:many ranges are equal-02} and \eqref{range contained in that of T-01}, it follows  that
\begin{equation}\label{2nd needed-01}U^*U\ell(A)=\ell(A), \quad U^*U(I_1-A)U^*=(I_1-A)U^*,\end{equation}
which leads to
\begin{equation*}(I_1-A)U^*U=U^*U(I_1-A)U^*U.\end{equation*}
Taking adjoints gives
\begin{equation}\label{equ:for exp of mq}(I_1-A)U^*U=U^*U(I_1-A).
\end{equation}
Since $A\ell(A)=\ell(A)A$, we have
\begin{equation*}-A\ell(A)U^*-\ell(A)(I_1-A)U^*=-\ell(A)U^*.
\end{equation*}
Based on the above equation, together with \eqref{equ:UpQUp*}, \eqref{2nd needed-01} and the second equation in \eqref{equ:RU*=R Q112-Q11 U*Q=0},  it is straightforward to verify that
\begin{equation*}(U_PQU_P^{*})^2=U_PQU_P^{*},
\end{equation*}
which implies that $Q^2=Q$,  so $Q$ is an idempotent.

From \eqref{equ:block matrix P},  we have
$$U_P(2P-I)U_P^*=\left(
                   \begin{array}{cc}
                     I_{1} & 0 \\
                     0 & -I_2 \\
                   \end{array}
                 \right).
$$
Combining this with \eqref{equ:UpQUp*}, which states  $U_PQ^*U_P^*=(U_PQU_P^*)^*$, we obtain
$$U_PQ^*U_P^*=U_P(2P-I)U_P^*\cdot U_PQU_P^*\cdot U_P(2P-I)U_P^*.$$
Consequently, $Q^*=(2P-I)Q(2P-I)$. Therefore, by Lemma~\ref{thm:short description of qpp}, $(P,Q)$ is a quasi-projection pair. Furthermore, it follows from \eqref{equ:new defn of T3 and T4}, \eqref{equ:block matrix P}, and \eqref{equ:UpQUp*}  that
\begin{equation}\label{equ:one range-1}U_P T_3 U_P^{*}=\left(
                    \begin{array}{cc}
                      0 & -\ell(A)U^* \\
                     0  & 0\\
                    \end{array}
                  \right).\end{equation}
Using \eqref{equ:many ranges are equal-02} and Lemma~\ref{lem:rang characterization-1},  we have
$$\overline{\mathcal{R}\big(-\ell(A)U^*\big)}=\overline{\mathcal{R}\big(\ell(A)U^*\big)}=\overline{\mathcal{R}[\ell^2(A)]}=\mathcal{R}(U^*U).$$
Applying \eqref{equ:one range-1}, this leads to
$$\overline{\mathcal{R}(U_PT_3U_P^{*})}=\mathcal{R}(U^*U)\oplus \{0\}.$$
As a result,
$$H_5=\overline{\mathcal{R}(T_3)}=\overline{\mathcal{R}(T_3U_P^*)}=U_P^{*}\,\overline{\mathcal{R}(U_PT_3U_P^{*})}=\mathcal{R}(U^*U).$$
Similarly, it can be shown that $H_6=\mathcal{R}(UU^*)$.
This shows that $H_5$ and $H_6$ are both orthogonally complemented in $H$, as
$\mathcal{R}(U^*U)$ and $\mathcal{R}(UU^*)$ are the ranges of  projections. Therefore,  according to Definition~\ref{defn of harmonious quasi-projection pair}, $(P,Q)$ is  harmonious.

$\Longrightarrow$. Suppose that $(P,Q)$ is a harmonious quasi-projection pair. Then $U_PQU_P ^*$ admits the block matrix decomposition \eqref{equ:block matrix Q}, characterized by    \begin{equation}\label{equ:Q11*=Q11 Q22*=Q22}
  A^*=A\in\mathcal{L}\big(\mathcal{R}(P)\big),\quad  Q_{22}^*=Q_{22}\in\mathcal{L}\big(\mathcal{N}(P)\big).\end{equation}
By \eqref{equ:existence of polar decomposition-01},  $Q_{21}\in\mathcal{L}\big(\mathcal{R}(P),\mathcal{N}(P)\big)$ admits  polar decomposition. Combining the idempotency condition $(U_PQU_P ^*)^2=U_PQU_P ^*$ with \eqref{equ:block matrix Q} and  \eqref{equ:Q11*=Q11 Q22*=Q22}, we derive the following operator identities:
\begin{align}
  \label{equ:Q21*Q21=A2-A} & Q_{21}^*Q_{21}=A^2-A, \\
   \label{equ:Q11Q21*+Q21*Q22=Q21*}&  AQ_{21}^*+Q_{21}^*Q_{22}=Q_{21}^*,\\
  \label{equ:-Q21Q21* +Q22^2=Q22} &    Q_{21}Q_{21}^*=Q_{22}^2-Q_{22}.
\end{align}
From \eqref{equ:Q21*Q21=A2-A}, the positivity $A^2-A\geq 0$ follows directly. Applying \eqref{defn of operator ell A} and \eqref{equ:Q21*Q21=A2-A}, we obtain $\ell(A)=|Q_{21}|$.
Lemmas~\ref{lem:Range closure of TT and T} and \ref{lem:Range Closure of T alpha and T} ensure:
$\overline{\mathcal{R}(Q_{21}^*)}=\overline{\mathcal{R}(|Q_{21}|)}$.
Thus, $Q_{21}$ admits  polar decomposition:
$$ Q_{21}=U\ell(A),\quad U^*U=P_{\overline{\mathcal{R}(Q_{21}^*)}}=P_{\overline{\mathcal{R}[\ell(A)]}},$$
implying $Q_{21}^*=\ell(A)U^*$,
$AQ_{21}^*=A\ell(A)U^*=\ell(A)AU^*$, and  the validity of \eqref{equ:many ranges are equal-02}. Substituting into
\eqref{equ:Q11Q21*+Q21*Q22=Q21*} yields
$$\ell(A)(U^*-AU^*-U^*Q_{22})=0.$$
Notably
\begin{align*}&\mathcal{R}(U^*Q_{22})\subseteq \mathcal{R}(U^*)=\overline{\mathcal{R}[\ell(A)]},
\end{align*}
and by  \eqref{range contained in that of T-01} we have  $$\overline{\mathcal{R}\big[(I_1-A)U^*\big]}=\overline{\mathcal{R}[\ell(A)]}.$$
Consequently, $$\mathcal{R}(U^*-AU^*-U^*Q_{22})\subseteq \overline{\mathcal{R}\big[(I_1-A)U^*\big]}+ \overline{\mathcal{R}(U^*Q_{22})}\subseteq \overline{\mathcal{R}[\ell(A)]}.$$
Hence
$$\mathcal{R}(U^*-AU^*-U^*Q_{22})\subseteq \overline{\mathcal{R}[\ell(A)]}\cap \mathcal{N}[\ell(A)]=\{0\}.$$
As a result,  $$U^*-AU^*-U^*Q_{22}=0,$$ which gives
\begin{equation}\label{equ:U(1-Q11)U*=UU*Q22}
 U(I_1-A)U^*=UU^*Q_{22}.
\end{equation}
Taking adjoints, we obtain
\begin{equation}\label{equ:UU*Q22=Q22UU*}
 UU^*Q_{22}=Q_{22}UU^*.
\end{equation}
Substituting $Q_{21}=U\ell(A)$ into \eqref{equ:-Q21Q21* +Q22^2=Q22} yields
\begin{equation}\label{equ:-U(Q11^2-Q11)U*+Q22^2=Q22}
  U(A^2-A)U^*=Q_{22}^2-Q_{22}.
\end{equation}
Define
\begin{equation}\label{equ:Q0=1-UU* Q22}
Q_0=(I_2-UU^*)Q_{22}.
\end{equation}
By construction, $U^*Q_{0}=0$. Using \eqref{equ:Q0=1-UU* Q22} and \eqref{equ:U(1-Q11)U*=UU*Q22}, we find
$$Q_{22}=UU^*Q_{22}+Q_0=U(I_1-A)U^*+Q_0.$$
From the second equation in \eqref{equ:Q11*=Q11 Q22*=Q22}, this means that
$$Q_0^*=\big[Q_{22}-U(I_1-A)U^*\big]^*=Q_0. $$
Additionally, from \eqref{equ:Q0=1-UU* Q22} and \eqref{equ:UU*Q22=Q22UU*}, we have
$$\quad Q_0^2=(I_2-UU^*)^2Q_{22}^2=(I_2-UU^*)Q_{22}^2.$$
Combining this with  \eqref{equ:-U(Q11^2-Q11)U*+Q22^2=Q22},  we conclude that
$$Q_0^2=(I_2-UU^*)\big[Q_{22}+U(A^2-A)U^*\big]=(I_2-UU^*)Q_{22}=Q_0.$$
This confirms that  $Q$  admits the canonical $2\times 2$ block matrix representation induced by $P$.
\end{proof}

\begin{remark}\label{remark of Q} Suppose that $(P,Q)$ is a harmonious quasi-projection pair with $P$ non-trivial. Let \eqref{equ:UpQUp*} denote the canonical $2\times 2$ block matrix representation of $Q$ induced by $P$. Define operators $\widetilde{U}\in \mathcal{L}\big(\mathcal{R}(P)\oplus \mathcal{R}(P),\mathcal{R}(P)\oplus\mathcal{N}(P)\big)$ and $\widetilde{Q_0}\in \mathcal{L}\big(\mathcal{R}(P)\oplus \mathcal{N}(P)\big)$ by
\begin{equation}\label{matrix widetitde U and Q0}
  \widetilde{U}=\left( \begin{array}{cc} I_{\mathcal{R}(P)} & 0 \\ 0 & U \\  \end{array}  \right),\quad \widetilde{Q_0}= \left( \begin{array}{cc} 0 & 0 \\ 0 & Q_0 \\  \end{array}  \right).
\end{equation}
Clearly, $\widetilde{U}$ is a partial isometry  and $\widetilde{Q_0}$ is a projection. With these,  $U_PQU_P^*$ can be rewritten  as
\begin{equation}\label{equ:rewrite UpQUp*}
  U_PQU_P^*=\widetilde{U}\widetilde{Q_1}\widetilde{U}^*+\widetilde{Q_0},
\end{equation}
where
\begin{equation}\label{equ:relate to rewite Q}
\widetilde{Q_1}=\left(
              \begin{array}{cc}
                A & -\ell(A) \\
                \ell(A) & I_{\mathcal{R}(P)}-A \\
              \end{array}
            \right)\in\mathcal{L}\big(\mathcal{R}(P)\oplus\mathcal{R}(P)\big).
\end{equation}
Let $U_0\in \mathcal{L}\big(\mathcal{R}(P)\oplus \mathcal{R}(P)\big)$ be defined by\footnote{For simplicity, we adopt the convention that for elements $a$ and $b_{ij}$ $(1\le i,j\le 2)$ in a commutative $C^*$-algebra, the multiplication of $a$  and the matrix $(b_{ij})$ is defined as
$$a\left(
                                                                                                                                \begin{array}{cc}
                                                                                                                                  b_{11} & b_{12} \\
                                                                                                                                  b_{21} & b_{22} \\
                                                                                                                                \end{array}
                                                                                                                              \right)=\left(
                                                                                                                                \begin{array}{cc}
                                                                                                                                  ab_{11} & ab_{12} \\
                                                                                                                                 a b_{21} & ab_{22} \\
                                                                                                                                \end{array}
                                                                                                                              \right).$$}
\begin{equation}\label{matrix unitary U0}
U_0=|2A-I_{\mathcal{R}(P)}|^{-\frac{1}{2}} \left( \begin{array}{cc} g(A) & -h(A) \\ -h(A) & -g(A) \\  \end{array}  \right),
\end{equation}
where the functions $g$ and $h$ are specified in \eqref{equ: the function of g,h,ell}. From the second equation in  \eqref{equ:newtemp-003}, $U_0$ is shown to be a unitary symmetry, i.e., a self-adjoint unitary operator. For every $t\in (0,1)^c$, the following identities hold:
\begin{equation*}
2t-1=f(t) |2t-1|, \quad h(t)g(t)=-f(t)\ell(t).
\end{equation*}
Therefore,
\begin{equation}\label{equ:|2A-1|=f(A)(2A-1)} 2A-I_{\mathcal{R}(P)}=f(A)|2A-I_{\mathcal{R}(P)}|,\quad  h(A)g(A)=-f(A)\ell(A).
\end{equation}
Combining \eqref{matrix unitary U0}, \eqref{equ:relate to rewite Q},  \eqref{equ:temp-002}, \eqref{equ:the properties of f and g}  and \eqref{equ:|2A-1|=f(A)(2A-1)}, we compute:
\begin{align*}
  U_0^*\widetilde{Q_1}U_0=&|2A-I_{\mathcal{R}(P)}|^{-\frac12}\left(
                                        \begin{array}{cc}
                                          g(A)(2A-I_{\mathcal{R}(P)}) & h(A)(2A-I_{\mathcal{R}(P)}) \\
                                          0 & 0 \\
                                        \end{array}
                                      \right)U_0\\
=&f(A)|2A-I_{\mathcal{R}(P)}|^{\frac12} \left(
                                        \begin{array}{cc}
                                          g(A) & h(A) \\
                                          0 & 0 \\
                                        \end{array}
                                      \right)U_0\\
=&f(A) \left(
                                        \begin{array}{cc}
                                          f(A) & 2f(A)\ell(A) \\
                                          0 & 0 \\
                                        \end{array}
                                      \right)=\left(
                                     \begin{array}{cc}
                                         I_{\mathcal{R}(P)} & 2\ell(A) \\
                                         0 & 0 \\
                                       \end{array}
                                    \right).
\end{align*}
Hence
\begin{equation}\label{equ:unitary equ to Q1}
  \widetilde{Q_1}=U_0\left(
                                        \begin{array}{cc}
                                          I_{\mathcal{R}(P)} & 2\ell(A) \\
                                          0 & 0 \\
                                        \end{array}
                                      \right)U_0^*.
\end{equation}
\end{remark}

Extending the canonical $2\times 2$ matrix representation established in Theorem~\ref{thm:2 by 2 canonical representation}, we prove the following characterization.

\begin{theorem}\label{thm:equivalence of unitary}Suppose that $(P,Q)$ is a harmonious quasi-projection with $P$ non-trivial. Let \eqref{equ:UpQUp*} denote the canonical $2\times 2$ block matrix representation of $Q$ induced by $P$. Then the following statements are  equivalent:
\begin{enumerate}
\item[{\rm (i)}] $\mathcal{N}(A)=\mathcal{N}(I_{\mathcal{R}(P)}-A)=\{0\}$ and $\mathcal{N}(U^*)=\{0\}$;
\item[{\rm (ii)}] $\mathcal{N}(A^2-A)=\{0\}$ and $\mathcal{N}(U^*)=\{0\}$;
\item[{\rm (iii)}] $U\in\mathcal{L}\big(\mathcal{R}(P),\mathcal{N}(P)\big)$ is a unitary operator;
\item[{\rm (iv)}] $\mathcal{N}(A)=\mathcal{N}(I_{\mathcal{R}(P)}-A)=\{0\}$ and $\mathcal{N}(Q_{22})=\mathcal{N}(I_{\mathcal{N}(P)}-Q_{22})=\{0\}$, where $$Q_{22}=U(I_{\mathcal{R}(P)}-A)U^*+Q_0.$$
\end{enumerate}
\end{theorem}
\begin{proof}(i)$\Longleftrightarrow$(ii). The commutativity of $A$ and $I_{\mathcal{R}(P)}-A$ implies  that both $\mathcal{N}(A)$ and $\mathcal{N}(I_{\mathcal{R}(P)}-A)$ are contained in $\mathcal{N}(A^2-A)$. Thus, the equivalence follows from the containment relations.

(ii)$\Longleftrightarrow$(iii). By \eqref{equ:RU*=R Q112-Q11 U*Q=0}, $\overline{\mathcal{R}(A^2-A)}$ is orthogonally complemented in $\mathcal{R}(P)$. Since $A$ is self-adjoint,  we have
$$\mathcal{R}(P)=\overline{\mathcal{R}(A^2-A)}+\mathcal{N}(A^2-A)=\mathcal{R}(U^*U)+ \mathcal{N}(A^2-A),$$
which means that $U^*U=I_{\mathcal{R}(P)}\Longleftrightarrow \mathcal{N}(A^2-A)=\{0\}$. Additionally, $UU^*=I_{\mathcal{N}(P)}\Longleftrightarrow \mathcal{N}(U^*)=\{0\}$.

(iii)+(i)$\Longrightarrow$(iv). It suffices to verify the second condition in (iv). Since $U$ is  unitary, \eqref{equ:RU*=R Q112-Q11 U*Q=0} implies $Q_0=0$. Therefore,
$$\mathcal{N}(Q_{22})=\mathcal{N}\big[U(I_{\mathcal{R}(P)}-A)U^*\big]=\{0\},\quad  \mathcal{N}(I_{\mathcal{N}(P)}-Q_{22})=\mathcal{N}(UAU^*)=\{0\}.$$

(iv)$\Longrightarrow$(i).  Suppose $x\in \mathcal{N}(U^*)$. Then $Q_{22}x=Q_0x$. Using  $U^*Q_0=0$ and $Q_0$ being a projection, we have
\begin{align*}(I_{\mathcal{N}(P)}-Q_{22})Q_{22}x=\big(I_{\mathcal{N}(P)}-Q_{0}-U(I_{\mathcal{R}(P)}-A)U^*\big)Q_0x=0,\end{align*}
which forces $x=0$ by the injectivity of $I_{\mathcal{N}(P)}-Q_{22}$ and $Q_{22}$.
\end{proof}

As a consequence of Theorem~\ref{thm:2 by 2 canonical representation},  we derive new $2\times 2$ block matrix representations for $P_{\mathcal{R}(Q)}$, $P_{\mathcal{N}(Q)}$, and $m(Q)$ in the remainder of this section.

\begin{corollary}\label{cor:derived rep for range Q}Suppose that $(P,Q)$ is a harmonious quasi-projection pair with $P$ non-trivial. Let \eqref{equ:UpQUp*} denote the canonical $2\times 2$ block matrix representation of $Q$ induced by $P$.  Then
\begin{equation}\label{equ:2nd form of range Q}U_PP_{\mathcal{R}(Q)} U_P^{*}=\left(
                    \begin{array}{cc}
                      B & B^\frac12\big(I_{\mathcal{R}(P)}-B\big)^\frac12 U_1^* \\
                     U_1B^\frac12\big(I_{\mathcal{R}(P)}-B\big)^\frac12  & U_1\big(I_{\mathcal{R}(P)}-B\big)U_1^*+Q_0\\
                    \end{array}
                  \right),
\end{equation}
where $B\in\mathcal{L}\big(\mathcal{R}(P)\big)$ is a positive contraction and $U_1\in \mathcal{L}\big(\mathcal{R}(P),\mathcal{N}(P)\big)$ is a partial isometry satisfying
\begin{equation}\label{two pre-conditions on the decomposition}\mathcal{R}(U_1^*)=\overline{\mathcal{R}\big(B(I_{\mathcal{R}(P)}-B)\big)}\quad\mbox{and}\quad U_1^*Q_0=0,
\end{equation}
in which the function $f$ is defined by  \eqref{defn of f-ist time}, and
\begin{equation}\label{B equals}
B=A\big(2A-I_{\mathcal{R}(P)}\big)^{-1},\quad U_1=Uf(A).
\end{equation}
\end{corollary}
\begin{proof}Put $I_{\mathcal{R}(P)}=I_1$,  $I_{\mathcal{N}(P)}=I_2$ and set
\begin{equation}\label{equ:Q21 and Q22}
Q_{21}=U\ell(A),\quad Q_{22}=U(I_1-A)U^*+Q_0.
\end{equation}
By \cite[Theorem~1.3]{Koliha},  we have
\begin{equation}\label{range projection Q}
P_{\mathcal{R}(Q)}=Q(Q+Q^*-I)^{-1}.
\end{equation}
Combining this  with \eqref{equ:block matrix Q} and \eqref{B equals} yields
\begin{equation*}U_P P_{\mathcal{R}(Q)}U_P^{*}=\left(
                                   \begin{array}{cc}
                                     B &  -Q_{21}^*(2Q_{22}-I_2)^{-1}\\
                                     Q_{21}(2A-I_1)^{-1} & (2Q_{22}-I_2)^{-1}Q_{22} \\
                                   \end{array}
                                 \right).
\end{equation*}
Since $U_P P_{\mathcal{R}(Q)}U_P^{*}$ is self-adjoint, this implies
\begin{equation}\label{equ:2nd pre form of range Q}U_P P_{\mathcal{R}(Q)}U_P^{*}=\left(
                                   \begin{array}{cc}
                                     B &  \big(Q_{21}(2A-I_1)^{-1}\big)^*\\
                                     Q_{21}(2A-I_1)^{-1} & (2Q_{22}-I_2)^{-1}Q_{22} \\
                                   \end{array}
                                 \right).
\end{equation}
Note that $f(A)$ is a unitary symmetry, so the operator $U_1$ defined by \eqref{B equals} is a partial isometry such that
$$U_1U_1^*=UU^*,\quad U_1^*Q_0=f(A)U^*Q_0=0.$$

Next, we verify the first condition in \eqref{two pre-conditions on the decomposition} and also show that
\begin{align}
  \label{equ:Q21(2A-I)inv} & Q_{21}(2A-I_1)^{-1}=U_1B^\frac12 (I_1-B)^\frac12, \\
  \label{equ:Q22(2Q22-I)inv}  & (2Q_{22}-I_2)^{-1}Q_{22}=U_1\big(I_1-B\big)U_1^*+Q_0.
\end{align}
For each $t\in\sigma(A)$, we have $t\in (0,1)^c$, so
\begin{align*}&0\le \frac{t}{2t-1}\le 1,\quad (2t-1)\sqrt{\frac{t}{2t-1}}\sqrt{1-\frac{t}{2t-1}}=f(t)\ell(t).\end{align*}
This confirms that $B$ is a positive contraction, and
\begin{equation}\label{needed for range-00}B^\frac12 (I_1-B)^\frac12(2A-I_1)=f(A)\ell(A),
\end{equation}
leading to
\begin{equation*}f(A)B^\frac12 (I_1-B)^\frac12=\ell(A)(2A-I_1)^{-1}.
\end{equation*}
The equation above, together with the first equation in \eqref{equ:Q21 and Q22}, yields
\begin{align*}
  & Q_{21}(2A-I_1)^{-1}=U\ell(A)(2A-I_1)^{-1}= U_1B^\frac12 (I_1-B)^\frac12.
  \end{align*}
Thus, \eqref{equ:Q21(2A-I)inv} is derived. Using \eqref{equ:many ranges are equal-02}, \eqref{needed for range-00} and the invertibility of $2A-I_1$, we have
 \begin{align*}\mathcal{R}(U_1^*)=&\overline{\mathcal{R}\big[f(A)U^*\big]}=\overline{\mathcal{R}\big[f(A)\ell(A)\big]}=\overline{\mathcal{R}\big[B^\frac12 (I_1-B)^\frac12(2A-I_1)\big]}\\
 =&\overline{\mathcal{R}\big[B^\frac12 (I_1-B)^\frac12\big]}=\overline{\mathcal{R}\big[B (I_1-B)\big]}.
\end{align*}
Therefore, the  first equation in \eqref{two pre-conditions on the decomposition} is established.

Finally, defining
 $$S=(2Q_{22}-I_2)\left[U_1\big(I_1-B\big)U_1^*+Q_0\right],$$ and utilizing $U_1=Uf(A)$, $f(A)B=Bf(A)$ and $f^2(A)=I_1$, we have $$U_1\big(I_1-B\big)U_1^*=U\big(I_1-B\big)U^*.$$
 Meanwhile, the second equation in \eqref{2nd needed-01} establishes $$U(I_1-A)U^*U=U(I_1-A).$$
 Combining with \eqref{equ:Q21 and Q22} and the identities
 $$U^*Q_0=0,\quad Q_0^2=Q_0,\quad (I_1-2A)(I_1-B)=I_1-A, $$  we verify
 \begin{align*}S=&(2Q_{22}-I_2)\left[U\big(I_1-B\big)U^*+Q_0\right]\\
 =&\left[2U(I_1-A)U^*+2Q_0-I_2\right]\left[U\big(I_1-B\big)U^*+Q_0\right]\\
 =&2U(I_1-A)\big(I_1-B\big)U^*-U(I_1-B)U^*+Q_0\\
=&U(I_1-2A)(I_1-B)U^*+Q_0\\
=&U(I_1-A)U^*+Q_0=Q_{22}.
\end{align*}
Hence, \eqref{equ:Q22(2Q22-I)inv} is valid. Moreover, \eqref{equ:2nd form of range Q} can be derived immediately by substituting \eqref{equ:Q21(2A-I)inv} and  \eqref{equ:Q22(2Q22-I)inv} into \eqref{equ:2nd pre form of range Q}.
\end{proof}

\begin{corollary}\label{cor:derived rep for kernel Q} Suppose that $(P,Q)$ is a harmonious quasi-projection pair with  $P$ non-trivial. Let \eqref{equ:UpQUp*} be the canonical $2\times 2$ block matrix representation of $Q$ induced by $P$. Then
\begin{equation}\label{equ:2nd form of kernel Q}U_PP_{\mathcal{N}(Q)}U_P^{*}=\left(
                    \begin{array}{cc}
                      I_{\mathcal{R}(P)}-B & B^\frac12\big(I_{\mathcal{R}(P)}-B\big)^\frac12 U_1^* \\
                     U_1B^\frac12\big(I_{\mathcal{R}(P)}-B\big)^\frac12  & U_1BU_1^*+Q_1\\
                    \end{array}
                  \right),
\end{equation}
where $B$ and $U_1$ are as  in Corollary~\ref{cor:derived rep for range Q},
and $Q_1=I_{\mathcal{N}(P)}-U_1U_1^*-Q_0$ is a projection satisfying $U_1^*Q_1=0$.
\end{corollary}
\begin{proof}We adopt the same notation as in the proof of Corollary~\ref{cor:derived rep for range Q}.
Replacing $Q$ with $I-Q$ in \eqref{range projection Q} gives
\begin{equation}\label{equ:relation of Q and P NQ}
P_{\mathcal{N}(Q)}=(Q-I)(Q+Q^*-I)^{-1}.
\end{equation}
Using \eqref{equ:relation of Q and P NQ}, \eqref{equ:block matrix Q},  the self-adjointness of $U_PP_{\mathcal{N}(Q)}U_P^{*}$, the first equation in \eqref{B equals} and \eqref{equ:Q21(2A-I)inv}, we derive
\begin{align*}
  U_PP_{\mathcal{N}(Q)}U_P^{*} &=\left(
                                   \begin{array}{cc}
                                     (A-I_1)(2A-I_1)^{-1} & \left[ Q_{21}(2A-I_1)^{-1}\right]^*\\
                                     Q_{21}(2A-I_1)^{-1} & (Q_{22}-I_2)(2Q_{22}-I_2)^{-1} \\
                                   \end{array}
                                 \right)  \\
   &=\left(
                                   \begin{array}{cc}
                                     I_1-B & B^\frac12\big(I_1-B\big)^\frac12 U_1^*\\
                                    U_1B^\frac12\big(I_1-B\big)^\frac12 & (Q_{22}-I_2)(2Q_{22}-I_2)^{-1} \\
                                   \end{array}
                                 \right).
\end{align*}
Since
$$(Q_{22}-I_2)(2Q_{22}-I_2)^{-1}=I_2-Q_{22}(2Q_{22}-I_2)^{-1},$$
it follows from \eqref{equ:Q22(2Q22-I)inv} that
\begin{align*}(Q_{22}-I_2)(2Q_{22}-I_2)^{-1}=&I_2-U_1(I_1-B)U_1^*-Q_0=U_1BU_1^*+Q_1.
\end{align*}
This confirms  the validity of \eqref{equ:2nd form of kernel Q}. Moreover, as $U_1$ is a partial isometry satisfying $U_1^*Q_0=0$, we have
$U_1^*Q_1=U_1^*(I_2-U_1U_1^*)-U_1^*Q_0=0.$  This completes the proof.
\end{proof}

To derive an analogous $2\times 2$ block matrix representation for the matched projection $m(Q)$, we need the following two lemmas.

\begin{lemma}\label{lem:m(uqu) equ umq u} Suppose that $Q\in \mathcal{L}(H)$ is an idempotent and $U\in \mathcal{L}(H)$ is a partial isometry. If $U^*UQ=QU^*U$, then  $m(UQU^*)=Um(Q)U^*$.
\end{lemma}
\begin{proof} Let $\mathfrak{A}$ denote the $C^*$-subalgebra of $\mathcal{L}(H)$ generated by $Q$.
Since $U^*UQ=QU^*U$, we have $Q^*U^*U=U^*UQ^*$ by taking adjoints. If follows that
$$U^*U a=aU^*U,\quad \forall a\in \mathfrak{A}.$$
Clearly,  $|Q|\in \mathfrak{A}$, and by \cite[Example~4.1]{Xu} we have $m(Q)\in \mathfrak{A}$. Hence,
\begin{equation*}U^*U |Q|=|Q|U^*U,\quad U^*Um(Q)=m(Q)U^*U.
\end{equation*}
By  the  above  argument, we conclude that $UQU^*$ is an idempotent, $Um(Q)U^*$ is a projection, and
$$|(UQU^*)|=U|Q|U^*.$$

Now, from \cite[Remark~3.10]{TXF}, we have
\begin{equation*}\label{equ:range of mq}\mathcal{R}\big(m(Q)\big)=\mathcal{R}(|Q|+Q),\quad \mathcal{R}\big(m(UQU^*)\big)=\mathcal{R}\big(|UQU^*|+UQU^*\big).\end{equation*}
Therefore,
\begin{align*}\mathcal{R}\big(m(UQU^*)\big)=&\mathcal{R}\big(U(|Q|+Q)U^*\big)=\mathcal{R}\big(U(|Q|+Q)U^*U\big)\\
=&\mathcal{R}\big(UU^*U(|Q|+Q)\big)=\mathcal{R}\big(U(|Q|+Q)\big)=\mathcal{R}\big(Um(Q)\big)\\
=&\mathcal{R}\big(UU^*Um(Q)\big)=\mathcal{R}\big(Um(Q)U^*U\big)=\mathcal{R}\big(Um(Q)U^*\big).
\end{align*}
This shows that the projections $m(UQU^*)$ and $Um(Q)U^*$ have identical ranges, so they are equal.
\end{proof}

\begin{lemma}\label{lem:m( P+ Q)=P+ m (Q)} Let $P\in\mathcal{L}(H)$ be a projection and $Q\in\mathcal{L}(H)$ be an idempotent with $PQ=QP=0$. Then $m(P+Q)=P+m(Q)$.
\end{lemma}
\begin{proof}A straightforward application of Lemma~\ref{lem:m(uqu) equ umq u} yields
$$\big(Pm(Q)\big)\big(Pm(Q)\big)^*=Pm(Q)P=m(PQP)=0,$$
which gives $Pm(Q)=0$. Hence, $m(Q)P=0$ and  $P+m(Q)$ is a projection such that
$$\mathcal{R}\big(P+m(Q)\big)=\mathcal{R}(P)+\mathcal{R}\big(m(Q)\big).$$
Moreover, since $PQ=QP=0$, we conclude that $P+Q$ is an idempotent, and
\begin{align*}&PQ^*=Q^*P=P|Q|=|Q|P=P|Q^*|=|Q^*|P=0,\\
& |P+Q|=P+|Q|,\quad |(P+Q)^*|=|P+Q^*|=P+|Q^*|.
\end{align*}
By \cite[Corollary~3.13]{TXF}, we have  $\mathcal{R}\big(m(Q)\big)=\mathcal{R}(T_Q)$, where
$$T_Q=|Q|+|Q^*|+Q+Q^*.$$
Since $T_Q$ is a self-adjoint operator with closed range, it follows from \cite[Lemma~5.7]{Luo-Moslehian-Xu} that $\mathcal{R}(T_Q^2)=\mathcal{R}(T_Q)$. Hence, for any $x,y\in H$, we have $T_Qy=T_Q^2 z$ for some $z\in H$, whence
\begin{align*}Px+T_Qy=(4P+T_Q)\left(\frac14 Px+T_Qz\right),
\end{align*}
which implies that
$$\mathcal{R}(P)+\mathcal{R}(T_Q)\subseteq \mathcal{R}(4P+T_Q)\subseteq \mathcal{R}(4P)+\mathcal{R}(T_Q)=\mathcal{R}(P)+\mathcal{R}(T_Q).$$
Consequently,
$$\mathcal{R}\big(P+m(Q)\big)=\mathcal{R}(P)+\mathcal{R}\big(m(Q)\big)=\mathcal{R}(P)+\mathcal{R}(T_Q)=\mathcal{R}(4P+T_Q).$$
Applying \cite[Corollary~3.13]{TXF} again, we obtain
\begin{align*}\mathcal{R}\big(m(P+Q)\big)=&\mathcal{R}\Big(|P+Q|+|(P+Q)^*|+P+Q+(P+Q)^*\Big)\\
=&\mathcal{R}(4P+T_Q)=\mathcal{R}\big(P+m(Q)\big).
\end{align*}
As shown in the proof of Lemma~\ref{lem:m(uqu) equ umq u}, equal ranges for these operators imply equality of the operators themselves. Thus, $m(P+Q)=P+m(Q)$.
\end{proof}

We now establish a new $2\times 2$ block matrix decomposition for the matched projection $m(Q)$.

\begin{theorem}\label{thm:matrix rep for m Q when that Q is given}Suppose that $(P,Q)$ is a harmonious quasi-projection pair with $P$ non-trivial. Let \eqref{equ:UpQUp*} be the canonical $2\times 2$ block matrix representation of $Q$ induced by $P$. Then
\begin{equation}\label{equ:representation of mQ}
  U_P\,m(Q)\,U_P^*=\left(
                 \begin{array}{cc}
                   \frac{1}{2}\big(I_{\mathcal{R}(P)}+f(A)\big) & 0 \\
                   0 & \frac{1}{2} U\big(I_{\mathcal{R}(P)}-f(A)\big)U^*+Q_0 \\
                 \end{array}
               \right),
\end{equation}
where the function $f$ is defined in \eqref{defn of f-ist time}.
\end{theorem}
\begin{proof}  We follow the notations from  Remark~\ref{remark of Q} and simplify $I_{\mathcal{R}(Q)}$ to $I_1$. By \eqref{equ:rewrite UpQUp*} and \eqref{equ:unitary equ to Q1}, we have
\begin{equation}\label{equ:for rep of mq+11}
  U_PQU_P^*=\widetilde{U}\widetilde{Q_1}\widetilde{U}^*+\widetilde{Q_0},\quad \widetilde{Q_1}=U_0\widetilde{Q_2}U_0^*,
\end{equation}
where
$$\widetilde{Q_2}=\left(
                                        \begin{array}{cc}
                                          I_{1} & 2\ell(A) \\
                                          0 & 0 \\
                                        \end{array}
                                      \right),$$
and $\widetilde{U}$ is a partial isometry, $\widetilde{Q_0}$ is a projection, $\widetilde{Q_1}$ is an idempotent, and $U_0$ is a unitary operator. These operators are defined by \eqref{matrix widetitde U and Q0}, \eqref{equ:relate to rewite Q}, and \eqref{matrix unitary U0}, respectively. From their definitions, combined with the second equation in \eqref{equ:RU*=R Q112-Q11 U*Q=0} and \eqref{equ:for exp of mq}, we derive
\begin{align*}\widetilde{U}^*\widetilde{U}\widetilde{Q_1}=\widetilde{Q_1}\widetilde{U}^*\widetilde{U}\quad\text{and}\quad \widetilde{U}\widetilde{Q_1}\widetilde{U}^*\widetilde{Q_0}=\widetilde{Q_0}\widetilde{U}\widetilde{Q_1}\widetilde{U}^*=0.
\end{align*}
It follows from \eqref{equ:for rep of mq+11} and Lemmas~\ref{lem:m(uqu) equ umq u}--\ref{lem:m( P+ Q)=P+ m (Q)} that
\begin{align}U_Pm(Q)U_P^*=&m(U_PQU_P^*)=m(\widetilde{U}\widetilde{Q_1}\widetilde{U}^*)+\widetilde{Q_0}=\widetilde{U}m(\widetilde{Q_1})\widetilde{U}^*+\widetilde{Q_0}\nonumber\\
\label{equ:for rep of mq+12}=&\widetilde{U}U_0\cdot m(\widetilde{Q_2})\cdot U_0^*\widetilde{U}^*+\widetilde{Q_0}.
\end{align}

By \eqref{defn of widetile Q111} and \eqref{defn of widetilde Q wrt A B}, we establish the matrix representation of $m(\widetilde{Q_2})$:
\begin{equation*}\label{defn of widetilde Q wrt A B--1}
m(\widetilde{Q_2})=\frac12 \left(\begin{array}{cc}
        (B+I_{1})B^{-1} &  2B^{-1}\ell(A)\\
        2\ell(A)B^{-1} &  4\ell(A)^2\big[B(B+I_{1})\big]^{-1} \\
      \end{array}\right),
\end{equation*}
where $$B=\big(4\ell^2(A)+I_1\big)^{\frac{1}{2}}.$$
Direct computation yields the following simplifications:
\begin{equation*}\label{equ:new B}
B=|2A-I_1|,\quad 4\ell(A)^2\big[B(B+I_{1})\big]^{-1}=(B-I_1)B^{-1}.
\end{equation*}
Following the notation convention from Remark~\ref{remark of Q}, we further simplify:
\begin{equation*}\label{defn of widetilde Q wrt A B--1}
m(\widetilde{Q_2})=\frac{1}{2}B^{-1} \left(
                                  \begin{array}{cc}
                                    B+I_1 & 2\ell(A) \\
                                    2\ell(A) & B-I_1 \\
                                  \end{array}
                                \right).
\end{equation*}
Applying \eqref{equ:the properties of f and g} and \eqref{equ:temp-002} with the identity $2A-I_1=f(A)B$, we compute:
\begin{align*}U_0m(\widetilde{Q_2})=&B^{-\frac12}\left( \begin{array}{cc} g(A) & -h(A) \\ -h(A) & -g(A) \\  \end{array}  \right)m(\widetilde{Q_2})=\frac12 B^{-\frac32}\left(
                                                                                                                                                                         \begin{array}{cc}
                                                                                                                                                                           T_{11} & T_{12} \\
                                                                                                                                                                           T_{21} & T_{22} \\
                                                                                                                                                                         \end{array}
                                                                                                                                                                       \right),
\end{align*}
where
\begin{align*}T_{11}=&g(A)(B+I_1)-2h(A)\ell(A)=g(A)(B+I_1)+g(A)(2A-2I_1)\\
=&g(A)(B+2A-I_1)=g(A)B\big(I_1+f(A)\big),\\
T_{12}=&2g(A)\ell(A)-h(A)(B-I_1)=-h(A)B\big(I_1+f(A)\big),\\
T_{21}=&-h(A)(B+I_1)-2g(A)\ell(A)=h(A)B\big(f(A)-I_1\big),\\
T_{22}=&-2h(A)\ell(A)-g(A)(B-I_1)=g(A)B\big(f(A)-I_1\big).
\end{align*}
Hence,
$$U_0m(\widetilde{Q_2})=
\left(
    \begin{array}{cc}
      \frac{1}{2}\big(I_1+f(A)\big) & 0 \\
      0 & \frac{1}{2}\big(I_1-f(A)\big) \\
    \end{array}
  \right)U_0,$$
which gives
\begin{equation*}
U_0m(\widetilde{Q_2})U_0^*=\left(
    \begin{array}{cc}
      \frac{1}{2}\big(I_1+f(A)\big) & 0 \\
      0 & \frac{1}{2}\big(I_1-f(A)\big) \\
    \end{array}
  \right).
\end{equation*}
Substituting this into \eqref{equ:for rep of mq+12} and recalling \eqref{matrix widetitde U and Q0} yields the desired expression \eqref{equ:representation of mQ} for $U_P\,m(Q)\,U_P^*$.
\end{proof}

The matched pairs demonstrate the following specific findings.
\begin{theorem}\label{thm:representation of Q realate to mQ} Let $Q\in\mathcal{L}(H)$ be a non-projection idempotent such that the matched pair $\big(m(Q),Q\big)$ is harmonious. Then
\begin{equation*}\label{equ:representation of Q realate to mQ}
  U_{m(Q)}QU_{m(Q)}^*=\left(
                 \begin{array}{cc}
                   A & -A^{\frac{1}{2}}\big(A-I_{H_1}\big)^{\frac{1}{2}}U^* \\
                   UA^{\frac{1}{2}}\big(A-I_{H_1}\big)^{\frac{1}{2}} & U(I_{H_1}-A)U^* \\
                 \end{array}
               \right),
\end{equation*}
where $H_1=\mathcal{R}\big(m(Q)\big)$, $A\in\mathcal{L}(H_1)$ is positive, and
$U\in \mathcal{L}(H_1,H_1^\bot)$ is a partial isometry satisfying $$A\geq I_{H_1}\quad\text{and}\quad  \mathcal{R}(U^*)=\overline{\mathcal{R}\big(A-I_{H_1}\big)}.$$
\end{theorem}
\begin{proof} For notational simplicity, denote $I_{H_1}$ by $I_1$. Since $Q$ is a non-projection idempotent,
the identity $$Q^*=\big(2m(Q)-I\big)Q\big(2m(Q)-I\big)$$
yields that $m(Q)$ is  non-trivial.  By Theorem~\ref{thm:2 by 2 canonical representation}, the operator $U_{m(Q)}QU_{m(Q)}^*$ admits the block matrix form  \eqref{equ:UpQUp*},
where $A\in\mathcal{L}(H_1)_{\mbox{sa}}$ satisfies $A^2-A\ge 0$, $\ell(A)$ is defined as in \eqref{defn of operator ell A},
$U\in \mathcal{L}\big(H_1,H_1^\bot\big)$ is a partial isometry, and $Q_0\in \mathcal{L}(H_1^\bot)$ is a projection fulfilling
\eqref{equ:RU*=R Q112-Q11 U*Q=0}. Applying Theorem~\ref{thm:matrix rep for m Q when that Q is given}, we conclude that
$U_{m(Q)}m(Q)U_{m(Q)}^*$ can be expressed as \eqref{equ:representation of mQ}.
On the other hand, it is clear that
$$U_{m(Q)}m(Q)U_{m(Q)}^*=\left(
                           \begin{array}{cc}
                             I_1 & 0 \\
                             0 & 0 \\
                           \end{array}
                         \right).$$
Equating the corresponding blocks, we obtain
$$\frac{1}{2}\big(I_1+f(A)\big)=I_1,\quad \frac{1}{2}U\big(I_1-f(A)\big)U^*+Q_0=0.$$
Therefore,
$$f(A)=I_1\quad\text{and}\quad Q_0=0.$$
Since $f$ is defined by \eqref{defn of f-ist time}, it follows from $f(A)=I_1$ that
$\sigma(A)\subseteq [1,+\infty)$.  Consequently,  $A\ge I_1$, $\ell(A)=A^{\frac{1}{2}}(A-I_1)^{\frac{1}{2}}$,
and
$$\mathcal{R}(A^2-A)=\mathcal{R}\big((A-I_1)A\big)=\mathcal{R}(A-I_1),$$
which simplifies \eqref{equ:RU*=R Q112-Q11 U*Q=0} to $\mathcal{R}(U^*)=\overline{\mathcal{R}(A-I_1)}$. This completes the proof.
\end{proof}

\section{The Halmos-like  block matrix representations for harmonious quasi-projection pairs}\label{sec:66 representations for qpp}
In this section, we investigate the Halmos-like   block  matrix representations for harmonious quasi-projection pairs.

\begin{lemma}\label{lem:property of harmonious quasi projection pair}{\rm\cite[Lemma~4.1 and Remark~4.2]{TXF-03}} Suppose that  $(P,Q)$ is a harmonious quasi-projection pair. Define the closed submodules
$H_i (1\le i\le 6)$ as in \eqref{eqn:defn of H1 and H4}--\eqref{eqn:defn of H5 and H6}.  Then each $H_i$ is orthogonally complemented in $H$, and the following decompositions hold:
\begin{equation}\label{equ:P=p1+p2+p5}
  P=P_{H_1}+P_{H_2}+P_{H_5},\quad  I-P=P_{H_3}+P_{H_4}+P_{H_6}.
\end{equation}
\end{lemma}

Let $(P,Q)$ be a harmonious quasi-projection pair on $H$. By \eqref{eqn:defn of H1 and H4}--\eqref{eqn:defn of H5 and H6} and Lemma~\ref{lem:another form of H1 and H4}, we have  $H_i\bot H_j$ for $i\ne j$.
For $1\le i\le 6$, denote the projection $P_{H_i}$ by $P_i$. Define the unitary operator $U_{P,Q}:H\rightarrow \bigoplus\limits_{i=1}^{6}H_{i}$ by
\begin{equation}\label{equ:defn of Upq}U_{P,Q}(x)=\Big(P_1(x),P_2(x),\cdots,P_6(x)\Big)^T, \quad x\in H.\end{equation}
It is straightforward to verify that  $$U_{P,Q}^*\Big((x_1,x_2,\cdots,x_6)^T\Big)=\sum_{i=1}^6 x_i,\quad x_i\in H_i, i=1,2,\cdots,6,$$
which implies that for every $T\in \mathcal{L}(H)$,
\begin{equation}\label{def:Upq T Upq*}
 U_{P,Q}TU_{P,Q}^*= \left(
    \begin{array}{cccc}
      P_1TP_1|_{H_1} & P_1TP_2|_{H_2} & \ldots & P_1TP_6|_{H_6} \\
      P_2TP_1|_{H_1} & P_2TP_2|_{H_2} & \ldots & P_2TP_6|_{H_6} \\
      \vdots & \vdots & \quad & \vdots \\
      P_6TP_1|_{H_1} & P_6TP_2|_{H_2} & \ldots & P_6TP_6|_{H_6} \\
    \end{array}
  \right).
\end{equation}

Our  first result of this section concerns a $6\times 6$ block matrix representation  for a general harmonious quasi-projection pair.
\begin{theorem}\label{thm:the 6x6 representation of p and q} Let $(P,Q)$ be a harmonious quasi-projection pair on $H$ with $P$ non-trivial. Define the closed submodules $H_i(1\leq i \leq 6)$ as in \eqref{eqn:defn of H1 and H4}--\eqref{eqn:defn of H5 and H6}, and let $U_{P,Q}$ be the unitary operator defined in \eqref{equ:defn of Upq}. Then $H_5\neq\{0\} $ if and only if $ H_6\neq\{0\} $, which occurs if and only if $Q\ne Q^*$. In this case, the following block matrix representations hold:
\begin{align}\label{decomposition2 of P}U_{P,Q} P U_{P,Q}^*&=I_{H_1}\oplus I_{H_2}\oplus 0_{H_3}\oplus 0_{H_4}\oplus I_{H_5}\oplus 0_{H_6},\\
\label{decomposition2 of Q}U_{P,Q} Q U_{P,Q}^*&=I_{H_1}\oplus 0_{H_2}\oplus I_{H_3}\oplus 0_{H_4}\oplus \widehat{Q},
\end{align}
where
\begin{equation}\label{decomposition of Q hat}
  \widehat{Q}=\left(
  \begin{array}{cc}
    A & -\ell(A)U^* \\
    U\ell(A) & U(I_{H_5}-A)U^* \\
  \end{array}
\right) \in \mathcal{L}(H_5\oplus H_6),
\end{equation}
in which $A\in \mathcal{L}(H_5)$ satisfies
\begin{equation}\label{equ:A2-A geq 0}
A^2-A\geq 0 \quad \text{and} \quad \overline{\mathcal{R}(A^2-A)}=H_5,
\end{equation}
$\ell(A)$ is defined by \eqref{defn of operator ell A}, and $U\in \mathcal{L}(H_5,H_6)$ is a unitary operator.
\end{theorem}
\begin{proof}For $1\leq i\leq 6$, denote $P_{H_i}$ and $I_{H_i}$ by $P_i$ and $I_i$, respectively. Let $T_3$ and $T_4$ be defined by \eqref{equ:new defn of T3 and T4}.
By Lemma~\ref{lem:H5 is 0}, we have
\begin{equation*}
  H_5 \ne \{0\} \Longleftrightarrow  H_6\ne \{0\}\Longleftrightarrow Q^*\ne Q.
\end{equation*}

Now, we assume $H_5\neq\{0\} $.  By \eqref{equ:P=p1+p2+p5}, the space $H$ admits an orthogonal decomposition $H=H_0+ H_0^\bot$, where
\begin{align}\label{defn of H 0}H_0=H_1+ H_2+ H_3+ H_4,\quad H_0^\bot=H_5+ H_6.
\end{align} From \eqref{equ:P=p1+p2+p5}, \eqref{eqn:defn of H1 and H4}, \eqref{eqn:defn of H2 and H3}, \eqref{alternative exps of H1 and H4} and \eqref{alternative exps of H2 and H3}, we deduce
\begin{equation}\label{equ:PPi QPi=0 or 1}
PP_i=\left\{
  \begin{array}{ll}
    P_i, & i=1, 2, 5,\\
   0, & i=3, 4, 6,
  \end{array}
\right.
\quad \text{and} \quad
QP_i=P_iQ= \left\{
  \begin{array}{ll}
    P_i, & i=1, 3, \\
    0, & i=2, 4.
  \end{array}
\right.
\end{equation}
Thus, the closed submodules  $H_0$ and $H_0^{\perp}$ are invariant under  $P$, and $H_0$ is also invariant under $Q$. It remains to show that  $H_0^{\perp}$ is  invariant under $Q$. Since $T_3$ is defined in \eqref{equ:new defn of T3 and T4}, we compute
$$PQT_3=PQ[I-(I-P)]Q(I-P)=T_3-T_3Q(I-P)=T_3(I-Q+QP),$$
and hence $$QT_3=PQT_3+(I-P)QT_3=T_3(I-Q+QP)+T_4Q(I-P).$$
By \eqref{eqn:defn of H5 and H6} and \eqref{defn of H 0}, this implies $QH_5\subseteq H_0^\bot$. Replacing $P$ by $I-P$, we have $QH_6 \subseteq H_0^{\perp}$. Therefore,
$ QH_0^{\perp}\subseteq H_0^{\perp}.$
This, together with \eqref{def:Upq T Upq*} and \eqref{equ:PPi QPi=0 or 1}, yields the decomposition \eqref{decomposition2 of Q},
where
\begin{equation}\label{equ:temp form of Q-01}
\widehat{Q}=\left(\begin{array}{cc}P_5QP_5|_{H_5}&P_5QP_6|_{H_6}\\P_6QP_5|_{H_5}&P_6QP_6|_{H_6}\end{array}\right)\in\mathcal{L}(H_5\oplus H_6).
\end{equation}
Moreover, a straightforward application of \eqref{def:Upq T Upq*} and \eqref{equ:PPi QPi=0 or 1} gives \eqref{decomposition2 of P}.

Finally, we verify that
$\widehat{Q}$ in \eqref{equ:temp form of Q-01} can be written as in  \eqref{decomposition of Q hat} and satisfies  \eqref{equ:A2-A geq 0}.
From \eqref{equ:P=p1+p2+p5} and the last two equations in \eqref{equ:PPi QPi=0 or 1}, we have
\begin{align}
  \label{equ:P5QP5=PQP-P1} & P_5QP_5=(P-P_1-P_2)Q(P-P_1-P_2)=PQP-P_1,\\
   \nonumber & P_5QP_6=(P-P_1-P_2)Q(I-P-P_3-P_4)=T_3.
\end{align}
Hence, $T_3P_6=T_3$ and $$\overline{\mathcal{R}(P_5QP_6|_{H_6})}=\overline{\mathcal{R}(T_3|_{H_6})}=\overline{\mathcal{R}(T_3)}.$$
The same argument shows that
$$\overline{\mathcal{R}(P_6QP_5|_{H_5})}=\overline{\mathcal{R}(T_4|_{H_5})}=\overline{\mathcal{R}(T_4)}.$$
Therefore, by the harmony of $(P,Q)$ and Lemma~\ref{lem:keep orthogonal complementity}, we conclude that
$\overline{\mathcal{R}(P_5QP_6|_{H_6})}$ and $\overline{\mathcal{R}(P_6QP_5|_{H_5})}$ are orthogonally complemented in $H_5$ and $H_6$, respectively.
Consequently, the idempotent $\widehat{Q}$  given by \eqref{equ:temp form of Q-01} has the form \eqref{equ:block matrix Q}, where
\begin{equation}\label{new defn of A and Q 22}A=P_5QP_5|_{H_5},\quad  Q_{21}=P_6QP_5|_{H_5},\quad Q_{22}=P_6QP_6|_{H_6}\end{equation}  such that $A\in\mathcal{L}(H_5)_{\mbox{sa}}$, $Q_{22}\in\mathcal{L}(H_6)_{\mbox{sa}}$, and $Q_{21}\in\mathcal{L}(H_5,H_6)$ has the polar decomposition (see \eqref{equ:existence of polar decomposition-01}).
Thus, by the proof of Theorem~\ref{thm:2 by 2 canonical representation}, $\widehat{Q}$ can be formulated as in \eqref{equ:UpQUp*} with $\mathcal{R}(P)$ and $\mathcal{N}(P)$  replaced by $H_5$ and $H_6$, respectively,  such that $A^2-A\ge 0$ and \eqref{equ:RU*=R Q112-Q11 U*Q=0} is satisfied.
Therefore, we only need to show that the associated partial isometry $U\in\mathcal{L}(H_5,H_6)$ is in fact a unitary operator, since in that case $\mathcal{R}(U^*)=H_5$, and $Q_0=0$ can be derived directly from $U^*Q_0=0$.

By (iii)$\Longleftrightarrow$(iv) in Theorem~\ref{thm:equivalence of unitary}, it suffices to prove that
\begin{align}\label{injectivity of A}&\mathcal{N}(A)=\mathcal{N}(I_5-A)=\{0\},\\
\label{injectivity of Q 22}&\mathcal{N}(Q_{22})=\mathcal{N}(I_6-Q_{22})=\{0\},
\end{align}
where $A$ and $Q_{22}$ are defined in \eqref{new defn of A and Q 22}.
Let $x\in \mathcal{N}(A)\subseteq H_5$. By \eqref{equ:P5QP5=PQP-P1} and the first equation in \eqref{new defn of A and Q 22}, we have
\begin{equation*}
0=Ax=(PQP-P_1)x=PQPx=PQx.
\end{equation*}
Since $QH_5\subseteq H_0^\bot$, it follows that
$$Qx\in \mathcal{N}(P)\cap \mathcal{R}(Q) \cap H_0^\bot =H_3\cap H_0^\bot=\{0\}.$$
Hence,
\begin{equation*}
x\in \mathcal{N}(Q)\cap H_5=\mathcal{N}(Q)\cap\mathcal{R}(P) \cap H_5=H_2\cap H_5=\{0\},
\end{equation*}
which shows $\mathcal{N}(A)=\{0\}$.

Now, let $x\in \mathcal{N}(I_5-A)\subseteq H_5$. From \eqref{equ:P5QP5=PQP-P1}, we get
\begin{equation*}
0=(I_5-A)x=x-PQx=P(I-Q)x.
\end{equation*}
Since $(I-Q)x=x-Qx\in H_0^\bot$, we obtain
$$(I-Q)x\in \mathcal{N}(P)\cap \mathcal{N}(Q)\cap H_0^\perp=H_4\cap H_0^\perp=\{0\},$$
which implies $x=Qx$. Therefore,
$$x\in \mathcal{R}(Q)\cap H_5 =\mathcal{R}(Q)\cap  \mathcal{R}(P)\cap H_5=H_1 \cap H_5=\{0\},$$
and thus $\mathcal{N}(I_5-A)=\{0\}$. This completes the proof of \eqref{injectivity of A}. The proof of
\eqref{injectivity of Q 22} is similar.
\end{proof}

In the following theorem, we present a new $4\times 4$ block matrix representation for a harmonious matched pair.

\begin{theorem}\label{thm:new representation for a general idempotent} Let $Q\in\mathcal{L}(H)$ be a non-projection idempotent such that the matched pair $\big(m(Q),Q\big)$ is harmonious. Let the closed submodules $H_i(i=1,4,5,6)$ and the unitary operator $U_{m(Q),Q}$ be defined by \eqref{eqn:defn of H1 and H4}, \eqref{eqn:defn of H5 and H6} and \eqref{equ:defn of Upq}, respectively, with $P$ replaced by $m(Q)$. Then
\begin{align}
  \label{equ:block matrix mQ and Q-mQ} & U_{m(Q),Q}\cdot m(Q)\cdot U_{m(Q), Q}^*=I_{H_1}\oplus 0_{H_4}\oplus I_{H_5}\oplus 0_{H_6}, \\
  \label{equ:block matrix mQ and Q-Q} & U_{m(Q),Q}\cdot Q\cdot U_{m(Q),Q}^*=I_{H_1}\oplus 0_{H_4}\oplus\widehat{Q},
\end{align}
where
\begin{equation}\label{decomposition of Q0 hat}
\widehat{Q}=\left(
                        \begin{array}{cc}
                          A & -A^{\frac{1}{2}}(A-I_{H_5})^{\frac{1}{2}}U^* \\
                          UA^{\frac{1}{2}}(A-I_{H_5})^{\frac{1}{2}} & U(I_{H_5}-A)U^* \\
                        \end{array}
                      \right)\in\mathcal{L}(H_5\oplus H_6),
\end{equation}
in which $U\in \mathcal{L}(H_5, H_6)$ is a unitary operator, and $A\in \mathcal{L}(H_5)$ satisfies
\begin{equation}\label{2conditions of A}A\geq I_{H_5},\quad \overline{\mathcal{R}(A-I_{H_5})}=H_5.\end{equation}
\end{theorem}
\begin{proof} By Lemma~\ref{lem:matched pair H2 and H3 zero}, we have  $H_2=\{0\}$ and $H_3=\{0\}$. Since $Q\ne Q^*$, Lemma~\ref{lem:H5 is 0} implies that $H_5\ne \{0\}$ and $H_6\ne \{0\}$. The desired conclusion follows directly from Theorem~\ref{thm:representation of Q realate to mQ} and the proof of Theorem~\ref{thm:the 6x6 representation of p and q}.
\end{proof}

We provide an additional result as follows.
\begin{theorem}\label{thm:new 6order rep for range and null}Let $Q\in\mathcal{L}(H)$ be a non-projection idempotent such that the matched pair $\big(m(Q),Q\big)$ is harmonious. Let the closed submodules $H_i(i=1,4,5,6)$ and the unitary operator $U_{m(Q),Q}$ be defined by \eqref{eqn:defn of H1 and H4}, \eqref{eqn:defn of H5 and H6} and \eqref{equ:defn of Upq}, respectively, with $P$ replaced by $m(Q)$. Then
\begin{align*}&U_{m(Q),Q}\cdot P_{\mathcal{R}(Q)}\cdot U_{m(Q),Q}^*=I_{H_1}\oplus  0_{H_4}\oplus \widehat{Q}_{\mathcal{R}},\\
&U_{m(Q),Q}\cdot P_{\mathcal{N}(Q)}\cdot U_{m(Q),Q}^*=0_{H_1}\oplus  I_{H_4}\oplus \widehat{Q}_{\mathcal{N}},
\end{align*}
where
\begin{align*}&\widehat{Q}_{\mathcal{R}}=\left(
                    \begin{array}{cc}
                      B & B^\frac12\big(I_{H_5}-B\big)^\frac12 U^* \\
                     U B^\frac12\big(I_{H_5}-B\big)^\frac12  & U\big(I_{H_5}-B\big)U^*\\
                    \end{array}
                  \right),\\
&\widehat{Q}_{\mathcal{N}}=\left(
                    \begin{array}{cc}
                      I_{H_5}-B & B^\frac12\big(I_{H_5}-B\big)^\frac12 U^* \\
                     U B^\frac12\big(I_{H_5}-B\big)^\frac12  & U B U^*\\
                    \end{array}
                  \right),
\end{align*}
in which $U\in \mathcal{L}(H_5, H_6)$ is a unitary operator, $A\in \mathcal{L}(H_5)$ satisfies \eqref{2conditions of A}, and $B=A(2A-I_{H_5})^{-1}\in \mathcal{L}(H_5)$ is a positive contraction.\end{theorem}
\begin{proof}We adopt the notations in Theorem~\ref{thm:new representation for a general idempotent} and Corollaries~\ref{cor:derived rep for range Q}--\ref{cor:derived rep for kernel Q}. Let the function $f$ be defined by \eqref{defn of f-ist time}. Since
$A\ge I_{H_5}$, we have $f(A)=I_{H_5}$, hence $U_1=Uf(A)=U$. Additionally, as $U_1$ is a unitary operator, we have $I_{H_6}-U_1U_1^*=0$.
The desired conclusion follows from \eqref{decomposition2 of Q}, \eqref{decomposition of Q hat}, \eqref{range projection Q}, \eqref{equ:relation of Q and P NQ}, and the proofs of Corollaries~\ref{cor:derived rep for range Q} and \ref{cor:derived rep for kernel Q}.
\end{proof}

\begin{remark}Suppose that  $A\in \mathcal{L}(H_5)$ satisfies \eqref{2conditions of A}. Let $B\in \mathcal{L}(H_5)$ be as in Theorem~\ref{thm:new 6order rep for range and null}. Since $B$ and $2A-I_{H_5}$ are invertible, we have
\begin{align*}\overline{\mathcal{R}\big(B(I_{H_5}-B)\big)}=&\overline{\mathcal{R}(I_{H_5}-B)}=\overline{\mathcal{R}\big((A-I_{H_5})(2A-I_{H_5})^{-1}\big)}\\
=&\overline{\mathcal{R}(A-I_{H_5})}=H_5.
\end{align*}
 \end{remark}

\section{Some applications}\label{sec:applications}
In this section, we present applications of Theorems~\ref{thm:the 6x6 representation of p and q} and \ref{thm:new representation for a general idempotent}.
\subsection{The supplementary projections and a new formula for an idempotent}\label{sec:representation for q}\label{subsec:supplementary projection}

Let $Q\in\mathcal{L}(H)$ be an idempotent. Since $\mathcal{L}(H)$ is a unital $C^*$-algebra, it can be identified with a $C^*$-subalgebra of $\mathbb{B}(K)$ for some Hilbert space $K$. Hence,  analogous to the Hilbert space setting, $Q$ can be represented using its range projection $P_{\mathcal{R}(Q)}$ and null space projection $P_{\mathcal{N}(Q)}$.  For example, the expression \eqref{formula for Q-001} is also valid in the Hilbert $C^*$-module case. Also,  invoking \cite[Theorem~3.8]{Ando02}, we have
$$Q=P_{\mathcal{R}(Q)}\big(P_{\mathcal{R}(Q)}+\lambda P_{\mathcal{N}(Q)}\big)^{-1}$$
for any $\lambda\ne 0$.  The primary objective of this subsection is to introduce a new projection $s(Q)$, termed the supplementary projection, and establish a new representation for every idempotent $Q$ via the two projections $m(Q)$ and $s(Q)$.

\begin{definition}\label{defn of sQ}Let $Q\in\mathcal{L}(H)$ be  an idempotent. The supplementary projection of $Q$, denoted by $s(Q)$, is defined as
\begin{equation}\label{equ:sQ is mQ}s(Q)=m\big(2P_{\mathcal{R}(Q)}-Q\big),
\end{equation}
where $m\big(2P_{\mathcal{R}(Q)}-Q\big)$ denotes the matched projection of the idempotent $2P_{\mathcal{R}(Q)}-Q$.
\end{definition}

For the remainder of this paper, the notation $s(Q)$ is reserved to denote the supplementary projection of an idempotent operator $Q$.
It is shown in \cite[Theorem~3.9]{TXF} that for every idempotent $Q\in\mathcal{L}(H)$,  $m(Q)$ is Murray-von Neumann equivalent to $P_{\mathcal{R}(Q)}$. The following theorem  establishes that $s(Q)$ admits the same equivalence relation.

\begin{theorem}\label{thm:range of sQ} For every idempotent $Q\in\mathcal{L}(H)$, there exist operators $T,V\in\mathcal{L}(H)$ such that
\begin{equation*}
  s(Q)=TT^\dag=VV^*\quad\mbox{and}\quad P_{\mathcal{R}(Q)}=T^\dag T=V^*V,
\end{equation*}
\end{theorem}
\begin{proof} Since $QP_{\mathcal{R}(Q)}=P_{\mathcal{R}(Q)}$, we have $P_{\mathcal{R}(Q)}Q^*=P_{\mathcal{R}(Q)}$ by taking adjoints. It follows that
\begin{align*}\big(2P_{\mathcal{R}(Q)}-Q\big)\big(2P_{\mathcal{R}(Q)}-Q\big)^*=QQ^*,
\end{align*}
which gives
\begin{equation}\label{same sqrate root}\big|\big(2P_{\mathcal{R}(Q)}-Q\big)^*\big|=|Q^*|.\end{equation}
Therefore, by \eqref{equ:final exp of P} we obtain
\begin{equation}\label{equ:def of sQ}
  s(Q)=\frac{1}{2}\big(|Q^*|+2P_{\mathcal{R}(Q)}-Q^*\big)|Q^*|^\dag (|Q^*|+I)^{-1} \big(|Q^*|+2P_{\mathcal{R}(Q)}-Q\big).
\end{equation}
As $Q$ is an idempotent, we have $\mathcal{R}(|Q^*|)=\mathcal{R}\big(|Q^*|^\dag\big)=\mathcal{R}(Q)$. Thus,
\begin{equation}\label{abs QQ dag}
P_{\mathcal{R}(Q)}=|Q^*|\cdot |Q^*|^\dag=|Q^*|^\dag \cdot |Q^*|,\quad Q|Q^*|=P_{\mathcal{R}(Q)}|Q^*|=|Q^*|.
\end{equation}
Taking adjoints gives $|Q^*|Q^*=|Q^*|P_{\mathcal{R}(Q)}=|Q^*|$.
Hence, $$|Q^*|^\dag (|Q^*|+I)=(|Q^*|+I)|Q^*|^\dag,$$ which implies
\begin{equation}\label{equ:commutativity of Qdag and Q+I}
  |Q^*|^\dag(|Q^*|+I)^{-1}=(|Q^*|+I)^{-1}|Q^*|^\dag.
\end{equation}
Combining this with \eqref{equ:def of sQ} yields
\begin{equation*}s(Q)=TW=VV^*,
\end{equation*}
where
\begin{align}\label{equ:new T with sQ}&T=|Q^*|+2P_{\mathcal{R}(Q)}-Q^*,\quad W=\frac{1}{2}|Q^*|^\dag(I+|Q^*|)^{-1}T^*,\\
&V=\frac{\sqrt{2}}{2}T\big(|Q^*|^\dag\big)^{\frac{1}{2}}\big(|Q^*|+I\big)^{-\frac{1}{2}}.\nonumber
\end{align}
Moreover,  we have
\begin{equation}\label{equ:pre wtt}
T^*T=2|Q^*|^2+2|Q^*|=2(I+|Q^*|)|Q^*|.
\end{equation}

Next, we prove $W=T^\dag$. Applying \eqref{equ:new T with sQ} and \eqref{equ:pre wtt} directly yields
\begin{equation*}WT=|Q^*|^\dag \cdot |Q^*|=P_{\mathcal{R}(Q)}.\end{equation*}
Observe that $\mathcal{R}(T^*)\subseteq \mathcal{R}(Q)$ and $\mathcal{R}(W)\subseteq \mathcal{R}(Q)$, hence
$$TWT=TP_{\mathcal{R}(Q)}=(P_{\mathcal{R}(Q)}T^*)^*=T,\quad WTW=P_{\mathcal{R}(Q)}W=W.$$
Consequently, we obtain
$$TWT=T,\quad WTW=W,\quad (TW)^*=TW,\quad (WT)^*=WT.$$
This establishes $W=T^\dag$, so
$$s(Q)=TT^\dag,\quad   P_{\mathcal{R}(Q)}=T^\dag T.$$

 Finally, using the definition of $V$ and \eqref{equ:pre wtt}, we compute
\begin{align*}
  V^*V & = (I+|Q^*|)^{-\frac{1}{2}}(|Q^*|^\dag)^{\frac{1}{2}}\cdot (I+|Q^*|)|Q^*| \cdot (|Q^*|^\dag)^{\frac{1}{2}} (I+|Q^*|)^{-\frac{1}{2}}\\
   & =|Q^*|^\dag \cdot |Q^*|=P_{\mathcal{R}(Q)}.
\end{align*}
This completes the proof.
\end{proof}

We present an alternative characterization of $s(Q)$ as follows.

\begin{theorem}\label{thm:relation of sQ and Q mQ} Let $Q\in\mathcal{L}(H)$ be an idempotent. Then
\begin{equation}\label{equ:equ:def of sQ-02}
  s(Q)=m(Q)+(I-Q^*)|Q^*|^\dag+|Q^*|^\dag(I-Q).
\end{equation}
\end{theorem}
\begin{proof} Observe that the operator $T$ defined by \eqref{equ:new T with sQ} can be rewritten  as
 $$T=|Q^*|+Q^*+2(P_{\mathcal{R}(Q)}-Q^*).$$
By \eqref{equ:commutativity of Qdag and Q+I}, \eqref{equ:final exp of P} and  \eqref{equ:def of sQ}, we have
\begin{equation*}\label{equ:median of sQ}X_1^*=X_1,\quad
s(Q)=m(Q)+2X_1+X_2+X_2^*,
\end{equation*}
where
\begin{align*}
 & X_1=(P_{\mathcal{R}(Q)}-Q^*)\cdot (|Q^*|+I)^{-1}|Q^*|^\dag\cdot (P_{\mathcal{R}(Q)}-Q),\\
   & X_2=(P_{\mathcal{R}(Q)}-Q^*)\cdot (|Q^*|+I)^{-1}|Q^*|^\dag \cdot (|Q^*|+Q).
 \end{align*}
A straightforward application of \eqref{abs QQ dag}  yields
\begin{align*}P_{\mathcal{R}(Q)}|Q^*|^\dag=|Q^*|^\dag=|Q^*|^\dag P_{\mathcal{R}(Q)}.
\end{align*}
Hence,
\begin{align*}|Q^*|^\dag\big(|Q^*|+P_{\mathcal{R}(Q)}\big)=|Q^*|^\dag\big(|Q^*|+I\big)=\big(|Q^*|+I\big)|Q^*|^\dag,
\end{align*}
and $(P_{\mathcal{R}(Q)}-Q^*)|Q^*|^\dag=(I-Q^*)|Q^*|^\dag$. Combining these results, we find
\begin{align*}
  X_1+X_2=(I-Q^*)|Q^*|^\dag.
\end{align*}
Since
$$2X_1+X_2+X_2^*=(X_1+X_2)+(X_1+X_2)^*,$$
substituting the above equality completes the proof.
\end{proof}

Building upon the derivation of the preceding two theorems, we further establish two propositions regarding the supplementary projections.

\begin{proposition}\label{prop:sq rep invariant}For every idempotent $Q\in\mathcal{L}(H)$, the following statements hold:
\begin{enumerate}
\item[{\rm (i)}] $s\big(2P_{\mathcal{R}(Q)}-Q\big)=m(Q)$;

\item[{\rm (ii)}]$\|s(Q)-m(Q)\|\le \|P_{\mathcal{R}(Q)}-Q\|$;

\item[{\rm (iii)}] If $X$ is a Hilbert space and $\pi:\mathcal{L}(H)\to\mathbb{B}(X)$ is a unital $C^*$-algebra morphism, then $s\big(\pi(Q)\big)=\pi\big(s(Q)\big)$.
\end{enumerate}
\end{proposition}
\begin{proof}(i) Let $Q_1=2P_{\mathcal{R}(Q)}-Q$. By \eqref{same sqrate root}, we have $|Q_1^*|=|Q^*|$, which implies
 $P_{\mathcal{R}(Q_1)}=P_{\mathcal{R}(Q)}$. It follows from \eqref{equ:sQ is mQ}  that
 $$s(Q_1)=m\big(2P_{\mathcal{R}(Q_1)}-Q_1\big)=m\big(2P_{\mathcal{R}(Q)}-Q_1\big)=m(Q).$$

 (ii) Let $Z=|Q^*|^\dag(I-Q)$. Then \eqref{equ:equ:def of sQ-02} can be rephrased as
 $$s(Q)-m(Q)=Z+Z^*.$$
Since $Z^2=0$, we have $Z^* Z\cdot ZZ^*=0$. Thus,
\begin{align*}\|s(Q)-m(Q)\|^2=&\|(Z+Z^*)^2\|=\|ZZ^*+Z^*Z\|\\
=&\max\{\|ZZ^*\|,\|Z^*Z\|\}=\|Z\|^2.
\end{align*}
Therefore,
\begin{align*}\|s(Q)-m(Q)\|=&\|Z\|=\big\| |Q^*|^\dag P_{\mathcal{R}(Q)}(I-Q)\big\|\\
\le&\big\|  |Q^*|^\dag\big\|\cdot \|P_{\mathcal{R}(Q)}-Q\|\le  \|P_{\mathcal{R}(Q)}-Q\|,
\end{align*}
where we have used \eqref{formula for the MP inverse-01} to conclude $\big\|  |Q^*|^\dag\big\|\le 1$.

(iii) Let $Q_\pi=\pi(Q)$. By \cite[Theorem~3.6]{TXF}, we have
$$m(Q_\pi)=\pi\big(m(Q)\big)\quad\mbox{and}\quad m\Big(\pi\big(2P_{\mathcal{R}(Q)}-Q\big)\Big)=\pi\Big(m\big(2P_{\mathcal{R}(Q)}-Q\big)\Big).$$
Due to \eqref{range projection Q}, we obtain $\pi(P_{\mathcal{R}(Q)})=P_{\mathcal{R}(Q_\pi)}$. It follows from \eqref{equ:sQ is mQ} that
\begin{align*}s(Q_\pi)=m\big(2P_{\mathcal{R}(Q_\pi)}-Q_\pi\big)=m\Big(\pi\big(2P_{\mathcal{R}(Q)}-Q\big)\Big)=\pi\big(s(Q)\big).
\end{align*}
This completes the proof.
\end{proof}

\begin{proposition} Let $Q\in\mathcal{L}(H)$ be an idempotent. Then $s(I-Q)=I-s(Q^*)$. Furthermore,
the following statements are equivalent:
\begin{enumerate}
  \item [\rm(i)] $s(I-Q)=I-s(Q)$;
  \item [\rm(ii)] $s(Q^*)=s(Q)$;
  \item [\rm(iii)] $Q$ is a projection;
  \item [\rm(iv)]$\big(s(Q),Q\big)$ is a quasi-projection pair.
\end{enumerate}
\end{proposition}
\begin{proof} By \eqref{equ:sQ is mQ} and  Lemma~\ref{thm:matched projection for I-Q}, we have
\begin{align*}s(Q^*)=&m\big(2P_{\mathcal{R}(Q^*)}-Q^*\big)=m\big(2P_{\mathcal{R}(Q^*)}-Q\big)=m\big(2(I-P_{\mathcal{N}(Q)})-Q\big)\\
=&m\big(I-\big[2P_{\mathcal{N}(Q)}-(I-Q)\big]\big)=I-m\big(2P_{\mathcal{N}(Q)}-(I-Q)\big)\\
=&I-m\big(2P_{\mathcal{R}(I-Q)}-(I-Q)\big)=I-s(I-Q).
\end{align*}
This establishes the identity $s(I-Q)=I-s(Q^*)$, thereby proving the equivalence of statements (i) and (ii).

The implication (iii)$\Longrightarrow$ (ii) is immediate.

(ii)$\Longrightarrow$(iii). We first note from Lemma~\ref{thm:matched projection for I-Q} that $m(Q^*)=m(Q)$. Substituting this into equation
\eqref{equ:equ:def of sQ-02} yields $X=Y$, where
$$X=(I-Q^*)|Q^*|^\dag+|Q^*|^\dag(I-Q),\quad Y=(I-Q)|Q|^\dag+|Q|^\dag(I-Q^*).$$
Consequently, we have
$$Q^*YQ=Q^*XQ=0.$$
Observe that
\begin{align*}Q^*(I-Q)|Q|^\dag Q=&\big(|Q|^\dag-|Q|^2 |Q|^\dag\big)Q=\big(|Q|^\dag-|Q|\big)Q=|Q|^\dag-|Q|,
\end{align*}
so \begin{align*}Q^*|Q|^\dag(I-Q^*) Q=&\big(Q^*(I-Q)|Q|^\dag Q\big)^*=\big(|Q|^\dag-|Q|\big)^*=|Q|^\dag-|Q|.
\end{align*}
Therefore,
$$2\big(|Q|^\dag-|Q|\big)=Q^*YQ=0,$$
and hence $$Q^*Q=|Q|^2=|Q|\cdot  |Q|^\dag=P_{\mathcal{R}(Q^*)}.$$
It follows that
$$\|Q\|^2=\big\|P_{\mathcal{R}(Q^*)}\big\|\le 1.$$
Since $Q$ is an idempotent, this norm condition holds if and only if $Q$ is a projection.

The implication (iii)$\Longrightarrow$(v) is straightforward, as $s(Q)=Q$ when $Q$ is a projection.

(iv)$\Longrightarrow$(iii). By Lemma~\ref{thm:short description of qpp} and \eqref{equ:sQ is mQ}, we have
\begin{align*}&Q^*=\big(2s(Q)-I\big)Q\big(2s(Q)-I\big),\\ &2P_{\mathcal{R}(Q)}-Q^*=\big(2s(Q)-I\big)\big(2P_{\mathcal{R}(Q)}-Q\big)\big(2s(Q)-I\big).\end{align*}
Hence,
$$P_{\mathcal{R}(Q)}\big(2s(Q)-I\big)=\big(2s(Q)-I\big)P_{\mathcal{R}(Q)},$$
which implies $P_{\mathcal{R}(Q)}s(Q)=s(Q)P_{\mathcal{R}(Q)}$. Equivalently,
\begin{equation}\label{self-adjoint01}\big(s(Q) |Q^*|\cdot |Q^*|^\dag\big)^*=s(Q) |Q^*|\cdot |Q^*|^\dag.\end{equation}
Meanwhile, by \eqref{equ:final exp of P} we have
$$m(Q)|Q^*|=\frac12(|Q^*|+Q^*).$$
Combining this equation with \eqref{equ:equ:def of sQ-02} yields
\begin{align*}s(Q) |Q^*|\cdot |Q^*|^\dag=&m(Q)|Q^*|\cdot |Q^*|^\dag+(I-Q^*)|Q^*|^\dag\\
=&\frac12 P_{\mathcal{R}(Q)}+|Q^*|^\dag-\frac12 Q^*|Q^*|^\dag.
\end{align*}
It follows from \eqref{self-adjoint01} that
$Q^*|Q^*|^\dag=|Q^*|^\dag Q$.
Consequently,
\begin{align*}Q=&|Q^*|\cdot |Q^*|^\dag Q=|Q^*|Q|Q^*|^\dag Q=|Q^*|QQ^*|Q^*|^\dag=|Q^*|^2,
\end{align*}
which implies $Q^*=Q$, so $Q$ is a projection.
\end{proof}

To derive the $4\times 4$ block matrix representations  for the supplementary projections, the following lemma is required.

\begin{lemma}\label{lem:rep of Qdag} Let $Q\in\mathcal{L}(H)$ be an idempotent. Then
\begin{equation}\label{equ:represetation of Qdag}|Q^*|=Q\big(2m(Q)-I\big),\quad
  |Q^*|^\dag=|Q^*|(Q+Q^*-I)^{-2}.
\end{equation}
\end{lemma}
\begin{proof} We begin by observing that
$$Q\big(|Q^*|+Q^*\big)=|Q^*|(I+|Q^*|).$$
Applying \eqref{equ:final exp of P}, we obtain
\begin{align*}Qm(Q)=\frac{1}{2}(|Q^*|+Q),
\end{align*}
which establishes the first identity in \eqref{equ:represetation of Qdag}. Consequently,
$$|Q^*|^\dag=\big(2m(Q)-I\big) Q^\dag,$$
since $2m(Q)-I$ is a unitary symmetry.

Next, we notice  that
$$(Q+Q^*-I)^{-1}Q=Q^*(Q+Q^*-I)^{-1}.$$
Therefore, by \eqref{range projection Q} we have
\begin{align*}Q^\dag=&Q^\dag Q\cdot QQ^\dag=P_{\mathcal{R}(Q^*)}P_{\mathcal{R}(Q)}=Q^*(Q+Q^*-I)^{-1}Q(Q+Q^*-I)^{-1}\\
=&Q^*(Q+Q^*-I)^{-2}.
\end{align*}
Consequently,
$$|Q^*|^\dag=\left[Q\big(2m(Q)-I\big)\right]^* (Q+Q^*-I)^{-2}=|Q^*|(Q+Q^*-I)^{-2}.$$
This completes the proof.
\end{proof}

The $4\times 4$ block matrix representation of a supplementary projection is given as follows.

\begin{theorem}\label{thm:rep4 sQ}Let $Q\in\mathcal{L}(H)$ be a non-projection idempotent such that the matched pair $\big(m(Q),Q\big)$ is harmonious. Let the closed submodules $H_i(i=1,4,5,6)$ and the unitary operator $U_{m(Q),Q}$ be defined by \eqref{eqn:defn of H1 and H4}, \eqref{eqn:defn of H5 and H6} and \eqref{equ:defn of Upq}, respectively, with $P$ replaced by $m(Q)$. Then
\begin{equation}\label{equ:6 times 6 of sq}
 U_{m(Q),Q}\cdot s(Q)\cdot U_{m(Q),Q}^*=I_{H_1}\oplus  0_{H_4}\oplus  \widehat{Q}_{s},
\end{equation}
where
\begin{equation}\label{equ:56 of sq}
  \widehat{Q}_{s}= \left( \begin{array}{cc}
   S & S^\frac12 (I_{H_5}-S)^\frac12 U^* \\
   US^\frac12 (I_{H_5}-S)^\frac12  & U(I_{H_5}-S)U^* \\
    \end{array}
   \right),
\end{equation}
in which $S=(2A-I_{H_5})^{-2}\in  \mathcal{L}(H_5)$.
\end{theorem}
\begin{proof} For notational simplicity, we denote $I_{H_i}=I_i$ for $(i=1,4,5,6)$. Let  $\widetilde{U}\in \mathcal{L}(H_5\oplus H_5, H_5\oplus H_6)$ be the unitary operator  defined by $$\widetilde{U}=\left(
                     \begin{array}{cc}
                       I_5 & 0 \\
                       0 & U \\
                     \end{array}
                   \right).
$$
By \eqref{equ:block matrix mQ and Q-mQ}--\eqref{decomposition of Q0 hat}, we have
\begin{align}
  \label{equ:6 times 6 of 1-qstar q} & U_{m(Q),Q}\cdot (I-Q^*)Q\cdot U_{m(Q),Q}^*=0\oplus 0 \oplus R_1,\\
  \label{equ:6 times 6 of 2mq-1} & U_{m(Q),Q}\cdot \big(2m(Q)-I\big)(Q+Q^*-I)^{-2}\cdot U_{m(Q),Q}^*=I_1\oplus I_4\oplus R_2.
\end{align}
where $A\ge I_5$ and
$$R_1=\widetilde{U}\left(
                                                       \begin{array}{cc}
                                                         -2A(A-I_5) & 2A^{\frac{1}{2}}(A-I_5)^{\frac{3}{2}}     \\
                                                         2A^{\frac{3}{2}}(A-I_5)^{\frac{1}{2}} & -2A(A-I_5) \\
                                                       \end{array}
                                                     \right)\widetilde{U}^*,\quad
 R_2=\widetilde{U}\left(
              \begin{array}{cc}
                S & 0 \\
                0 & -S \\
              \end{array}
            \right)\widetilde{U}^*.  $$
Following the convention established in Remark~\ref{remark of Q}, we compute:
$$R_1R_2=\widetilde{U}\cdot 2S  \left(
                                                       \begin{array}{cc}
                                                         -A (A-I_5)  & -A^{\frac{1}{2}}(A-I_5)^{\frac{3}{2}} \\
                                                         A^{\frac{3}{2}}(A-I_5)^{\frac{1}{2}} & A (A-I_5)  \\
                                                       \end{array}
                                                     \right)\cdot \widetilde{U}^*,$$
and hence,
\begin{equation*}\label{equ:R1R2 add R1R2star}
R_1R_2+(R_1R_2)^*=\widetilde{U}\cdot 2S   \left(
                                                       \begin{array}{cc}
                                                         -2A (A-I_5)  & A^{\frac{1}{2}}(A-I_5)^{\frac{1}{2}} \\
                                                         A^{\frac{1}{2}}(A-I_5)^{\frac{1}{2}} & 2A (A-I_5)  \\
                                                       \end{array}
                                                     \right)\cdot \widetilde{U}^*.
\end{equation*}
Since
\begin{align*}&I_5-4SA(A-I_5)=S\big(S^{-1}-4A(A-I_5)\big)=S,\\
&2SA^{\frac{1}{2}}(A-I_5)^{\frac{1}{2}}=S^\frac12 (I_5-S)^\frac12,\quad 4SA (A-I_5)=I_5-S,
\end{align*}
we have
\begin{align}
 \label{equ:56 of sq-02}  &  (I_5\oplus 0)+R_1R_2+(R_1R_2)^*  =\widehat{Q}_s,
\end{align}
where $\widehat{Q}_s$ is defined by \eqref{equ:56 of sq}.

Now, we  define
\begin{equation}\label{equ:simplic 1-q qdag}
 X=(I-Q^*)|Q^*|^\dag.
\end{equation}
By \eqref{equ:represetation of Qdag}, this becomes
$$X=(I-Q^*)\cdot Q\big(2m(Q)-I\big)\cdot (Q+Q^*-I)^{-2}.$$
Combining with \eqref{equ:6 times 6 of 1-qstar q} and \eqref{equ:6 times 6 of 2mq-1}, we obtain
$$U_{m(Q),Q}\cdot X\cdot U_{m(Q),Q}^*=0\oplus 0\oplus R_1R_2.$$
Hence,
\begin{equation}\label{equ:6 times 6 of X+Xstar}
U_{m(Q),Q}\cdot(X+X^*)\cdot U_{m(Q),Q}^*=0\oplus 0\oplus \big[R_1R_2+(R_1R_2)^*\big].
\end{equation}
Therefore, \eqref{equ:6 times 6 of sq} can be derived immediately from \eqref{equ:equ:def of sQ-02}, \eqref{equ:block matrix mQ and Q-mQ}, \eqref{equ:simplic 1-q qdag}, \eqref{equ:6 times 6 of X+Xstar} and \eqref{equ:56 of sq-02}.
\end{proof}

Based on Theorem~\ref{thm:rep4 sQ}, we now provide a new formula for an idempotent in terms of its matched projection and supplementary projection.

\begin{theorem}\label{thm:new formula for q bansed on mQ and sQ} Let $Q\in \mathcal{L}(H)$ be an idempotent. Then
\begin{equation}\label{equ:representation of Q}
  Q=\frac{1}{2} C^{-1}\Big[\big(C^{\frac{1}{2}}+s(Q)\big)\cdot \big(2m(Q)-I\big)+I-m(Q)\Big],
\end{equation}
where
\begin{equation}\label{equ:rep of C}
C=\big(I-m(Q)-s(Q)\big)^2.
\end{equation}
\end{theorem}
\begin{proof} If $Q$ is a projection, then $m(Q)=s(Q)=Q$, so $C=I$. Hence, \eqref{equ:representation of Q} holds trivially. Suppose now that $Q$ is a non-projection idempotent.

First, we consider the special case where the matched pair $\big(m(Q),Q\big)$ is harmonious.
From Equations \eqref{equ:block matrix mQ and Q-mQ} and \eqref{equ:6 times 6 of sq}, it follows that
$$U_{m(Q),Q}\cdot\big(I-m(Q)-s(Q)\big)\cdot U_{m(Q),Q}^*=-I_1\oplus I_4\oplus R_3,$$
where $$R_3=\widetilde{U}\cdot  \left(
                                           \begin{array}{cc}
                                         -S & -S^{\frac12}(I_5-S)^{\frac12} \\
                                      -S^{\frac12}(I_5-S)^{\frac12} & S \\
                                       \end{array}
                                    \right)
\cdot \widetilde{U}^*.$$
Hence,
$$U_{m(Q),Q}\cdot C\cdot U_{m(Q),Q}^*=I_1\oplus I_4\oplus \widetilde{U}\left(
                                                                         \begin{array}{cc}
                                                                          S &  \\
                                                                            & S \\
                                                                         \end{array}
                                                                       \right)\widetilde{U}^*.
$$
This implies that $C$ is a positive definite operator on $H$. Therefore,
\begin{align*}
  U_{m(Q),Q}\cdot C^{-\frac{1}{2}}\big(2m(Q)-I\big)\cdot U_{m(Q),Q}^*  = &I_1\oplus -I_4\oplus \widetilde{U}\left(
                                                                         \begin{array}{cc}
                                                                          S^{\frac12} &  \\
                                                                            & -S^{\frac12} \\
                                                                         \end{array}
                                                                       \right)\widetilde{U}^*.
  \end{align*}
On the other hand, since $S^{-\frac12}=2A-I_5$, by \eqref{equ:block matrix mQ and Q-mQ}--\eqref{decomposition of Q0 hat}, we have
\begin{align*}U_{m(Q),Q}\cdot (Q+Q^*-I)\cdot U_{m(Q),Q}^*=&I_1\oplus -I_4 \oplus \widetilde{U}
  \left(
    \begin{array}{cc}
      S^{-\frac12} & 0 \\
      0 & -S^{-\frac12} \\
    \end{array}
  \right)\widetilde{U}^*.
\end{align*}
This establishes that
\begin{equation}\label{equ:relation of Q and C -01}
 Q+Q^*= C^{-\frac{1}{2}}\big(2m(Q)-I\big)+I.
\end{equation}
Similarly, by Equations  \eqref{equ:block matrix mQ and Q-mQ}--\eqref{decomposition of Q0 hat} and \eqref{equ:6 times 6 of sq}--\eqref{equ:56 of sq}, one can also obtain
\begin{equation}\label{equ:relation of Q and C -02}
  Q-Q^*=C^{-1}\big(s(Q)m(Q)-m(Q)s(Q)\big).
\end{equation}
Moreover, since $s(Q)$ and $m(Q)$ are projections, by \eqref{equ:rep of C} it is easy to get $$C+s(Q)m(Q)-m(Q)s(Q)=s(Q)\big(2m(Q)-I\big)+I-m(Q),$$
which is combined with \eqref{equ:relation of Q and C -01} and \eqref{equ:relation of Q and C -02} to obtain
\begin{align*}
  Q & =\frac{1}{2}\Big[C^{-\frac{1}{2}}\big(2m(Q)-I\big)+I+C^{-1}\Big(s(Q)m(Q)-m(Q)s(Q)\Big)\Big] \\
   & =\frac{1}{2} C^{-1}\Big[C^{\frac{1}{2}}\big(2m(Q)-I\big)+C+s(Q)m(Q)-m(Q)s(Q)\Big]\\
   &=\frac{1}{2} C^{-1}\Big[C^{\frac{1}{2}}\big(2m(Q)-I\big)+s(Q)\big(2m(Q)-I\big)+I-m(Q)\Big]\\
   &=\frac{1}{2} C^{-1}\Big[\big(C^{\frac{1}{2}}+s(Q)\big)\cdot \big(2m(Q)-I\big)+I-m(Q)\Big].
\end{align*}
This confirms the validity of Equation \eqref{equ:representation of Q}.

Next, we consider the general case where the matched pair $\big(m(Q),Q\big)$ is not harmonious.
Since $\mathcal{L}(H)$ is a unital $C^*$-algebra,  we may fix a Hilbert space $X$ and a faithful unital $C^*$-algebra homomorphism $\pi:\mathcal{L}(H)\to \mathbb{B}(X)$. We use $I$ to denote the identity operator on both $H$ and $X$. Let $Q_\pi=\pi(Q)$. Then $\big(m(Q_\pi),Q_\pi\big)$ is a harmonious matched pair. By \cite[Theorem~3.6]{TXF} and
Proposition~\ref{prop:sq rep invariant}, we have
$$m(Q_\pi)=\pi\big(m(Q)\big), \quad s(Q_\pi)=\pi\big(s(Q)\big).$$ Thus,
\begin{equation*}
  Q_\pi=\frac{1}{2} C_\pi^{-1}\Big[\Big(C_\pi^{\frac{1}{2}}+\pi\big(s(Q)\big)\Big)\cdot \Big(2\pi\big(m(Q)\big)-I\Big)+I-\pi\big(m(Q)\big)\Big],
\end{equation*}
where $C_\pi=\big[I-\pi\big(m(Q)\big)-\pi\big(s(Q)\big)\big]^2$.
The desired identity \eqref{equ:representation of Q} now follows from the faithfulness of $\pi$.
\end{proof}

\subsection{A new proof of the canonical form for quadratic operators}

Recall that an operator $T\in\mathbb{B}(H)$ is called quadratic if it satisfies
\begin{equation}\label{defn of qot}(T-a I)(T-b I)=0\end{equation} for some complex numbers $a$ and $b$. Notable examples of quadratic operators include idempotents and square-zero operators. A key structural result, established in \cite[Theorem~1.1]{TW}, states that every quadratic operator in $\mathbb{B}(H)$ is unitarily equivalent to a block matrix of the form
$$aI\oplus b I\oplus \left(
                       \begin{array}{cc}
                         aI & B \\
                         0& bI\\
                       \end{array}
                     \right),$$
where $B$ is a positive operator with trivial kernel. The purpose of this subsection is to provide a new proof of this characterization.
\begin{lemma}\label{lem:unitary form of idem}Let $Q\in\mathbb{B}(H)$ be an idempotent. Then there exist Hilbert spaces $K_1,K_2,K_3$ and a unitary operator
$W: H\to K_1\oplus K_2\oplus K_3\oplus K_3$ such that
\begin{equation}\label{block for idempotent}WQW^*=I_{K_1}\oplus 0_{K_2}\oplus \left(
                                                    \begin{array}{cc}
                                                      I_{K_3} & B \\
                                                      0_{K_3} & 0_{K_3} \\
                                                    \end{array}
                                                  \right),\end{equation}
where $B\in\mathbb{B}(K_3)$ is a positive operator with $\mathcal{N}(B)=\{0\}$.
\end{lemma}
\begin{proof}First, suppose $Q$ is a projection. Define $K_1=\mathcal{R}(Q)$, $K_2=\mathcal{N}(Q)$, $K_3=\{0\}$ and $B=0$. For any $x\in H$, set
$$W(x)=\big(Qx,(I-Q)x,0,0\big)^T.$$ Clearly, $W: H\to K_1\oplus K_2\oplus K_3\oplus K_3$ is a unitary operator such that  \eqref{block for idempotent} is satisfied.

Next, consider the case where $Q\in\mathbb{B}(H)$ is a non-projection idempotent. Since $H$ is a Hilbert space, it follows that $\big(m(Q),Q\big)$ is harmonious. So by \eqref{equ:block matrix mQ and Q-Q}--\eqref{2conditions of A}, \eqref{equ:relate to rewite Q}  and  \eqref{equ:unitary equ to Q1}, we have
\begin{equation}\label{4norm Q01}U_{m(Q),Q}\cdot Q\cdot U_{m(Q),Q}^*=I_{H_1}\oplus 0_{H_4}\oplus\widehat{Q},\end{equation}
where
\begin{equation}\label{4norm Q02}\widehat{Q}=\left(
                              \begin{array}{cc}
                                I_{H_5} &  \\
                                 & U \\
                              \end{array}
                            \right)U_0\left(
                                        \begin{array}{cc}
                                          I_{H_5} & 2\ell(A) \\
                                          0 & 0 \\
                                        \end{array}
                                      \right)U_0^*\left(
                              \begin{array}{cc}
                                I_{H_5} &  \\
                                 & U^* \\
                              \end{array}
                            \right),
\end{equation}
in which $U\in \mathbb{B}(H_5, H_6)$ is a unitary operator, $A\ge I_{H_5}$, $\mathcal{N}(A-I_{H_5})=\{0\}$, $\ell(A)$ is defined by \eqref{defn of operator ell A} and $U_0\in\mathbb{B}(H_5\oplus H_5)$ is given as in \eqref{matrix unitary U0}. Define
$$K_1=H_1,\quad K_2=H_4,\quad K_3=H_5,\quad W=I_{H_1}\oplus I_{H_4}\oplus \left(
                                                                            \begin{array}{cc}
                                                                              I_{H_5} &  \\
                                                                               & U \\
                                                                            \end{array}
                                                                          \right)U_0.$$
Then $W: H\to K_1\oplus K_2\oplus K_3\oplus K_3$ is a unitary operator and satisfies  \eqref{block for idempotent} with $B=2\ell(A)$.
Since $B=2A^\frac12 (A-I_{H_5})^\frac12$, we have $B\ge 0$, and
$\mathcal{N}(B)=\mathcal{N}(A-I_{H_5})=\{0\}$. This completes the proof.
\end{proof}

\begin{lemma}\label{lem:square-zerot}Let $S\in\mathbb{B}(H)$ be a square-zero operator. Then there exist Hilbert spaces $K_1,K_2,K_3$ and a unitary operator
$W: H\to K_1\oplus K_2\oplus K_3\oplus K_3$ such that
\begin{equation*}\label{block for square-zerot}WSW^*=0_{K_1}\oplus 0_{K_2}\oplus \left(
                                                    \begin{array}{cc}
                                                      0_{K_3} & B \\
                                                      0_{K_3} & 0_{K_3} \\
                                                    \end{array}
                                                  \right),\end{equation*}
where $B\in\mathbb{B}(K_3)$ is a positive operator with $\mathcal{N}(B)=\{0\}$.
\end{lemma}
\begin{proof}If $S$ is self-adjoint, then $S^*S=S^2=0$, so $S=0$. The conclusion holds trivially in this case.

Now suppose $S$ is not self-adjoint. Let $P$ denote the projection from $H$ onto $\overline{\mathcal{R}(S)}$, and set $Q=P+S$. Then $Q\in\mathbb{B}(H)$ is a non-projection idempotent satisfying
\begin{equation}\label{rela of pq}QP=P,\quad PQ=Q. \end{equation}
By Lemma~\ref{lem:unitary form of idem}, there exist Hilbert spaces $K_1,K_2,K_3$ and a unitary operator
$W: H\to K_1\oplus K_2\oplus K_3\oplus K_3$ such that \eqref{block for idempotent} is satisfied. Write
$WPW^*$ as a block matrix $(P_{ij})_{1\le i,j\le 4}.$
Combining this with $P=P^*$, \eqref{block for idempotent} and \eqref{rela of pq}, a direct computation shows
$$WPW^*=I_{K_1}\oplus 0_{K_2}\oplus I_{K_3}\oplus 0_{K_3}.$$
Since $WSW^*=WQW^*-WPW^*$, the desired conclusion follows.
\end{proof}

\begin{theorem}\label{thm:quadratic form}{\rm \cite[Theorem~1.1]{TW}} Let $T\in\mathbb{B}(H)$ be a quadratic operator satisfying \eqref{defn of qot}. Then there exist Hilbert spaces $K_1,K_2,K_3$ and a unitary operator
$W: H\to K_1\oplus K_2\oplus K_3\oplus K_3$ such that
\begin{equation*}WTW^*=a I_{K_1}\oplus bI_{K_2}\oplus \left(
                                                    \begin{array}{cc}
                                                      a I_{K_3} & B \\
                                                      0_{K_3} & bI_{K_3} \\
                                                    \end{array}
                                                  \right),\end{equation*}
where $B\in\mathbb{B}(K_3)$ is a positive operator with $\mathcal{N}(B)=\{0\}$.
\end{theorem}
\begin{proof}If $a=b$, set $S=T-aI$, then $S$ is a square-zero operator. If $a\ne b$, define
$$Q=\frac{1}{a-b}(T-bI),$$ then $Q$ is an idempotent in $\mathbb{B}(H)$. Let $a-b=|a-b|e^{i\theta}$ for some $\theta\in [0,2\pi)$. For any Hilbert space $K$ and $B\in\mathbb{B}(K)$, we have
$$\left(
    \begin{array}{cc}
      I_K & 0 \\
      0 & e^{i\theta}I_K \\
    \end{array}
  \right)\left(
           \begin{array}{cc}
             I_K & (a-b)B \\
             0 & 0 \\
           \end{array}
         \right)\left(
    \begin{array}{cc}
      I_K & 0 \\
      0 & e^{i\theta}I_K \\
    \end{array}
  \right)^*=\left(
              \begin{array}{cc}
                I_K & |a-b|B \\
                0 & 0 \\
              \end{array}
            \right).$$
The desired conclusion follows by combining this matrix identity with Lemmas~\ref{lem:unitary form of idem}  and \ref{lem:square-zerot}.
\end{proof}

\subsection{Non-validity of the Krein-Krasnoselskii-Milman equality for general quasi-projection pairs}

For every projections $P,Q\in\mathcal{L}(H)$, by \cite[Lemma~4.1]{XY2} we have
\begin{equation*}\label{equ:norm of subtraction 2 projections}\|P-Q\|=\max\big\{\|P(I-Q)\|, \|(I-P)Q\| \big\},
\end{equation*}
which is known as the Krein-Krasnoselskii-Milman equality (or KKM-equality for short). It is interesting to investigate the validity of
this equality in the context of quasi-projection pairs. To this end, we require the following lemma.

\begin{lemma}\label{lem:technique-02 for norm} Let $A\in \mathcal{L}(H)$ be self-adjoint and satisfy $A^2-A\geq 0$. Let $\ell(A)$ denote the operator defined by \eqref{defn of operator ell A}, and let $\widetilde{Q_1}\in\mathcal{L}(H\oplus H)$ be defined as in
\eqref{equ:relate to rewite Q} (with $\mathcal{R}(P)$ replaced by $H$). Then
\begin{equation}\label{equ: another norm of Q and 1-Q}
\big\|I\oplus 0- \widetilde{Q_1}\big\|=\big\||A-I|+\ell(A)\big\|,\quad \Big\|\big(I\oplus 0- \widetilde{Q_1}\big)^2\Big\|=\|A-I\|.
\end{equation}
\end{lemma}
\begin{proof}Define $\Omega=I\oplus 0- \widetilde{Q_1}$. Then
\begin{align*}&\Omega=\left(
                                                 \begin{array}{cc}
                                                   I-A & \ell(A) \\
                                                   -\ell(A) & A-I \\
                                                 \end{array}
                                               \right),\quad \Omega^2=\left(
                                                                        \begin{array}{cc}
                                                                          I-A & 0 \\
                                                                          0 & I-A \\
                                                                        \end{array}
                                                                      \right),\\
&U_0^*\Omega U_0=\left(
                    \begin{array}{cc}
                      0 & I-A+\ell(A) \\
                      I-A-\ell(A) & 0 \\
                    \end{array}
                  \right),
\end{align*}
where $U_0\in\mathcal{L}(H\oplus H)$ is a unitary operator defined by
\begin{equation*}U_0=\left(
        \begin{array}{cc}
          \frac{\sqrt{2}}{2}I & \frac{\sqrt{2}}{2}I \\
          -\frac{\sqrt{2}}{2}I & \frac{\sqrt{2}}{2}I\\
        \end{array}
      \right).
\end{equation*}
Hence, $\|\Omega^2\|=\|A-I\|$, and
\begin{equation*}\|\Omega\|=\|U_0^*\Omega U_0\|=\max\big\{\|I-A+\ell(A)\|,\|A-I+\ell(A)\|\big\}.
\end{equation*}

For  $t\in\sigma(A)$, define
$$x(t)=|1-t+\ell(t)|,\quad y(t)=|t-1+\ell(t)|.$$
Then
\begin{align*}\|\Omega\|=&\max\big\{\max_{t\in\sigma(A)}x(t),\max_{s\in\sigma(A)}y(s)\big\}=\max_{t\in\sigma(A)}\max\{x(t),y(t)\}\\
=&\max_{t\in\sigma(A)}\big(|t-1|+\ell(t)\big)=\big\||A-I|+\ell(A)\big\|.
\end{align*}
This completes the proof.
\end{proof}

The next theorem shows that the KKM-equality  does not hold in general for quasi-projection pairs.

\begin{theorem}\label{thm:modified KKM formula}Let $(P,Q)$ be a quasi-projection pair on $H$. Then
\begin{align}\label{2pairs norm}&\|(I-P)Q\|=\|Q(I-P)\|, \quad \|(I-Q)P\|=\|P(I-Q)\|,\\
&\label{2sides inequality}\big\|(P-Q)^2\big\|\le \max\big\{\|(I-P)Q\|, \|(I-Q)P\|\big\}\le \|P-Q\|.
\end{align}
\end{theorem}
\begin{proof}Observe  that $2P-I$ is a unitary symmetry. By Lemmas~\ref{thm:four equivalences} and \ref{thm:short description of qpp}, we have
\begin{align*}\|P(I-Q)\|=\big\|\big(P(I-Q)\big)^*\big\|=\|(2P-I)(I-Q)(2P-I)P\|=\|(I-Q)P\|.
\end{align*}
Replacing $(P,Q)$ by $(I-P,I-Q)$ yields $\|(I-P)Q\|=\|Q(I-P)\|$.

Next, we prove  \eqref{2sides inequality}. If $Q$ is a projection, then \eqref{2sides inequality}
reduces to the KKM-equality. Thus, \eqref{2sides inequality} holds trivially.

Suppose now that $Q$ is a non-projection idempotent. Up to a faithful unital representation of $\mathcal{L}(H)$, we may assume $(P,Q)$ is harmonious.
Following the notation in the proof of Theorem~\ref{thm:the 6x6 representation of p and q}, we have
$H_5\ne \{0\}$ and $H_6\ne \{0\}$ by Lemma~\ref{lem:H5 is 0}, and
\begin{align*}&U_{P,Q} (P-Q)U_{P,Q} ^*=0\oplus I_2\oplus (-I_3)\oplus 0\oplus S_1,\\
&U_{P,Q} (I-P)QU_{P,Q} ^*=0\oplus 0\oplus I_3\oplus 0\oplus S_2,\\
&U_{P,Q} (I-Q)PU_{P,Q} ^*=0\oplus I_2\oplus 0\oplus 0\oplus S_3,
\end{align*}
where
\begin{align*}&S_1=\left(
                                                 \begin{array}{cc}
                                                   I_5-A & \ell(A)U^* \\
                                                   -U\ell(A) & -U(I_5-A)U^* \\
                                                 \end{array}
                                               \right)=\widetilde{U}\Omega\widetilde{U}^*,\\
&S_2=\widetilde{U}\left(
       \begin{array}{cc}
         0 & 0 \\
         \ell(A) & I_5-A \\
       \end{array}
     \right)\widetilde{U}^*,\quad S_3=\widetilde{U}\left(
       \begin{array}{cc}
         I_5-A & 0 \\
         -\ell(A) & 0 \\
       \end{array}
     \right)\widetilde{U}^*,
\end{align*}
in which $\ell(A)$ and $\widetilde{Q_1}$  are defined in \eqref{defn of operator ell A} and \eqref{equ:relate to rewite Q} with $\mathcal{R}(P)$ and $I$ replaced by $H_5$ and $I_5$, respectively, and
\begin{equation*}\widetilde{U}=\left(
                \begin{array}{cc}
                  I_5 & 0 \\
                  0 & U \\
                \end{array}
              \right), \quad
\Omega=I_5\oplus 0-\widetilde{Q_1}=\left(
                                                 \begin{array}{cc}
                                                   I_5-A & \ell(A) \\
                                                   -\ell(A) & A-I_5 \\
                                                 \end{array}
                                               \right).
\end{equation*}
Define
\begin{equation}\label{defn of M}M=\big\|(A-I_5)^2+\ell^2(A)\big\|^\frac12.\end{equation}
From \eqref{equ: another norm of Q and 1-Q},  we get $\big\|S_1^2\big\|=\big\|\Omega^2\big\|=\|A-I_5\|$, and
\begin{align*}\|S_1\|=&\|\Omega\|=\big\||A-I_5|+\ell(A)\big\|\ge M=\|S_2\|=\|S_3\|\ge \big\|S_1^2\big\|.
\end{align*}
Thus,
\begin{align*}\big\|(P-Q)^2\big\|=&\max\left\{\|I_2\|,\|I_3\|, \big\|S_1^2\big\|\right\}\le  \max\left\{\|I_2\|,\|I_3\|, M\right\}\\
=&\max\big\{\max\{\|I_3\|, M\}, \max\{\|I_2\|, M\}\big\}\\
=&\max\big\{\|(I-P)Q\|, \|(I-Q)P\|\big\}\\
\le&\max\left\{\|I_2\|,\|I_3\|, \|S_1\|\right\}=\|P-Q\|.
\end{align*}
This confirms the validity of \eqref{2sides inequality}.
\end{proof}

For the matched pairs, the preceding theorem can be specified as follows.
\begin{theorem}\label{thm:KKME for matched pair}Let  $Q\in\mathcal{L}(H)$ be a non-projection idempotent. Then the operators
$$\big(I-m(Q)\big)Q,\quad Q\big (I-m(Q)\big),\quad (I-Q)m(Q),\quad m(Q)(I-Q)$$
share the same norm, equal to $\frac{\sqrt{2}}{2}\sqrt{\|Q\|(\|Q\|-1)}$. Furthermore,
\begin{align}&\Big\|\big(m(Q)-Q\big)^2\Big\|=\frac{\|Q\|-1}{2},\nonumber\\
\label{mq qdistance}&\big\|m(Q)-Q\big\|=\frac12\Big(\|Q\|-1+\sqrt{\|Q\|^2-1}\Big).
\end{align}
Consequently, the following strict inequalities hold:
$$\sqrt{2}\Big\|\big(m(Q)-Q\big)^2\Big\|<\big\|\big(I-m(Q)\big)Q\big\|<\big\|m(Q)-Q\big\|<\sqrt{2}\big\|\big(I-m(Q)\big)Q\big\|.$$
\end{theorem}
\begin{proof}By a faithful unital representation of $\mathcal{L}(H)$, we may assume $\big(m(Q),Q\big)$ is harmonious. We use the notation from the proof of Theorem~\ref{thm:modified KKM formula}. From \eqref{2conditions of A},  $A\geq I_{H_5}$. By \eqref{4norm Q01} and \eqref{4norm Q02},  $\|Q\|=\max\{\|I_{H_1}\|,\|\widehat{Q}\|\}$, and
$$\|\widehat{Q}\|=\left\|\left(
                                \begin{array}{cc}
                                  I_{H_5} & 2\ell(A) \\
                                  0 & 0 \\
                                \end{array}
                              \right)\right\|=2\|A\|-1.$$
Since $Q$ is not a projection, $\|Q\|>1$, so $\|Q\|=\|\widehat{Q}\|$. This gives
\begin{equation}\label{norm of a wrt q}\|A\|=\frac{\|Q\|+1}{2}>1.\end{equation}
The number $M$ defined in \eqref{defn of M} thus simplifies to
$$M=\sqrt{(\|A\|-1)(2\|A\|-1)}.$$
By Lemma~\ref{lem:matched pair H2 and H3 zero}, $H_2=\{0\}$ and $H_3=\{0\}$. Invoking the proof of Theorem~\ref{thm:modified KKM formula}, we find
\begin{align*}&\big\|\big(I-m(Q)\big)Q\big\|=\big\|(I-Q)m(Q)\big\|=M, \\
&\big\|m(Q)-Q\big\|=\|A\|-1+\sqrt{\|A\|(\|A\|-1)},\\
&\Big\|\big(m(Q)-Q\big)^2\Big\|=\|A\|-1.
\end{align*}
From $\|A\|>1$, it is straightforward to verify
$$\sqrt{2}(\|A\|-1)<M<\|A\|-1+\sqrt{\|A\|(\|A\|-1)}<\sqrt{2}M.$$
The desired conclusion follows from \eqref{2pairs norm} and \eqref{norm of a wrt q}.
\end{proof}

\begin{remark}Formula \eqref{mq qdistance} was originally derived in \cite[Theorem~4.17]{TXF} using a different approach. For an alternative proof of this formula, see also  \cite[Theorem~4.4]{ZFL}.
\end{remark}

\vspace{5ex}










\begin{thebibliography}{99}



\bibitem{Afriat} S. Afriat, Orthogonal and oblique projectors and the characteristics of pairs of vector spaces, Proc. Cambridge Philos. Soc. 53 (1957), 800--816.


\bibitem{Ando02} T. Ando, Unbounded or bounded idempotent operators in Hilbert space, Linear Algebra Appl 438 (2013), no. 10, 3769--3775.
















\bibitem{BS}A. B\"ottcher and I. M. Spitkovsky, A gentle guide to the basics of two projections theory, Linear Algebra Appl. 432 (2010), no. 6, 1412--1459.


\bibitem{Corach-Maestripieri}G. Corach and A. Maestripieri, Products of orthogonal projections and polar decompositions, Linear Algebra Appl. 434 (2011), no. 6, 1594--1609.









\bibitem{FXY}C. Fu, Q. Xu and G. Yan, Characterizations of the harmonious pairs of projections on a Hilbert $C^*$-module, Linear Multilinear Algebra 70 (2022), no. 21, 6312--6320.








\bibitem{Halmos} P. R. Halmos, Two subspaces, Trans. Amer. Math. Soc. 144 (1969), 381--389.










\bibitem{Koliha} J. J. Koliha, Range projections of idempotents in $C^*$-algebras, Demonstratio Math. 34 (2001), no. 1, 91--103.



\bibitem{Lance}E. C. Lance, Hilbert $C^*$-modules--A toolkit for operator algebraists, Cambridge University Press, Cambridge, 1995.






\bibitem{Liu-Luo-Xu}N. Liu, W. Luo and Q. Xu, The polar decomposition for adjointable operators on Hilbert $C^*$-modules and centered operators, Adv. Oper. Theory 3 (2018), no. 4, 855--867.



\bibitem{Luo-Moslehian-Xu}W. Luo, M. S. Moslehian and Q. Xu, Halmos' two projections theorem for Hilbert $C^*$-module operators and the Friedrichs angle of two closed submodules, Linear Algebra Appl. 577 (2019), 134--158.



\bibitem{MT} V. M. Manuilov and E. V. Troitsky, Hilbert $C^*$-modules, Translated from the 2001 Russian original by the authors, Translations of Mathematical Monographs, 226, American Mathematical Society, Providence, RI, 2005.







\bibitem{Ovchinnikov}P. G. Ovchinnikov, Automorphisms of the poset of skew projections, J. Funct. Anal. 115 (1993), no. 1, 184--189.

\bibitem{Paschke}W. L. Paschke, Inner product modules over $B^*$-algebras, Trans. Amer. Math. Soc. 182 (1973), 443--468.













\bibitem{TXF}X. Tian, Q. Xu and C. Fu, The matched projections of idempotents on Hilbert $C^*$-modules, J. Operator Theory 94 (2025), no. 1, 101--133.
 doi: 10.7900/jot.2023oct21.2450



\bibitem{TXF-03}X. Tian, Q. Xu and C. Fu, Characterizations of the semi-harmonious and harmonious quasi-projection pairs on Hilbert $C^*$-modules, Quaest. Math. (2025), 1--24.
 doi:10.2989/16073606.2025.2550026

\bibitem{TW} S. Tso and P. Wu, Matricial ranges of quadratic operators, Rocky Mountain J.
Math. 29 (1999), 1139--1152.

\bibitem{Vosough-Moslehian-Xu} M. Vosough, M. S. Moslehian and Q. Xu,
Closed range and nonclosed range adjointable operators on Hilbert $C^*$-modules,
Positivity 22 (2018), 701--710.




\bibitem{Xu}Q. Xu, A new formula for the weighted Moore-Penrose inverse and its applications, Banach J. Math. Anal. 19 (2025), no. 4, Paper No. 55.

\bibitem{XS}Q. Xu and L. Sheng, Positive semi-definite matrices of adjointable operators on Hilbert $C^*$-modules,
Linear Algebra Appl. 428 (2008), no. 4, 992--1000.




\bibitem{Xu-Wei-Gu}Q. Xu, Y. Wei and Y. Gu, Sharp norm estimations for Moore-Penrose inverses of stable
perturbations of Hilbert $C^*$-module operators, SIAM J. Numer. Anal. 47 (2010), 4735--4758.

\bibitem{Xu-Yan}Q. Xu and G. Yan, Harmonious projections and Halmos' two projections theorem for Hilbert $C^*$-module operators, Linear Algebra Appl. 601 (2020), 265--284.

\bibitem{XY2} Q. Xu and G. Yan, Products of projections, polar decompositions and norms of differences of two projections, Bull. Iranian Math. Soc. 48 (2022), no. 1, 279--293.




\bibitem{ZFL}C. Zhao, Y. Fang and Y. Li, Some characterizations of the quasi-projection pairs and the matched projections, Banach J. Math. Anal. 19 (2025), no. 2, Paper No. 18.


\end{thebibliography}
\end{document}